\documentclass[11pt, fleqn]{amsart}

\usepackage[square,compress,comma, numbers,sort]{natbib}
\usepackage[colorlinks=true, citecolor=red, linkcolor=blue]{hyperref}
\usepackage{amsfonts,mathtools}
\usepackage{amsmath}

\allowdisplaybreaks[4]

\usepackage{amssymb}
\usepackage{color}
\usepackage{float}

\DeclarePairedDelimiterXPP\pk[1]{\mathbb{P}}\{ \}{}{ #1}
\DeclarePairedDelimiterXPP\E[1]{\mathbb{E}}\{ \}{}{	#1}

\newcounter{alphasect}
\def\alphainsection{0}

\let\oldsection=\section
\def\section{%
  \ifnum\alphainsection=1%
    \addtocounter{alphasect}{1}
  \fi%
\oldsection}%

\renewcommand\thesection{%
 \ifnum\alphainsection=1%
   \Alph{alphasect}%
 \else
   \arabic{section}%
 \fi%
}%

\newenvironment{alphasection}{%
  \ifnum\alphainsection=1%
    \errhelp={Let other blocks end at the beginning of the next block.}
    \errmessage{Nested Alpha section not allowed}
  \fi%
  \setcounter{alphasect}{0}
  \def\alphainsection{1}
}{%
  \setcounter{alphasect}{0}
  \def\alphainsection{0}
}%

\usepackage{xparse}% http://ctan.org/pkg/xparse

\NewDocumentCommand{\ceil}{s O{} m}{%
	\IfBooleanTF{#1} % starred
	{\left\lceil#3\right\rceil} % \ceil*[..]{..}
	{#2\lceil#3#2\rceil} % \ceil[..]{..}
}
\NewDocumentCommand{\floor}{s O{} m}{%
	\IfBooleanTF{#1} % starred
	{\left\lfloor#3\right\rfloor}
	{#2\lfloor#3#2\rfloor}
}

\def\x{\vk{x}}

\definecolor{c20}{rgb}{0.,0.7,0.}
\definecolor{c30}{rgb}{0.,0.,1.}
\definecolor{c40}{rgb}{1,0.1,0.7}
\definecolor{c50}{rgb}{1,0,0}
\definecolor{c60}{rgb}{1,0.9,0.1}
\definecolor{c70}{rgb}{0.50,1.00,0.00}

\numberwithin{equation}{section}
\newtheorem{theo}{Theorem}[section]
\newtheorem{sat}[theo]{Proposition}
\newtheorem{de}[theo]{Definition}
\newtheorem{lem}[theo]{Lemma}

\newtheorem{example}[theo]{Example}
\newtheorem{korr}[theo]{Corollary}
\newtheorem{remark}[theo]{Remark}

\numberwithin{equation}{section}

\newcommand{\prooftheo}[1]{ \textsc{Proof of Theorem} \ref{#1} }

\newcommand{\QED}{\hfill $\Box$}

\newcommand{\COM}[1]{}

\def\IF{\infty}

\newcommand{\R}{\mathbb{R}}
\newcommand{\inr}{\in \R}

\newcommand{\BQN}{\begin{eqnarray}}
\newcommand{\EQN}{\end{eqnarray}}
\newcommand{\BQNY}{\begin{eqnarray*}}
	\newcommand{\EQNY}{\end{eqnarray*}}

\newcommand{\limit}[1]{\lim_{#1 \to   \infty}}

\newcommand{\kb}[1]{\boldsymbol{#1}}
\newcommand{\vk}[1]{\kb{#1}}

\def\bqny#1{{\begin{eqnarray*} #1 \end{eqnarray*}}}
\def\bqn#1{{\begin{eqnarray} #1 \end{eqnarray}}}

\newcommand{\BS}{\begin{sat}}
	\newcommand{\ES}{\end{sat}}
\newcommand{\BT}{\begin{theo}}
	\newcommand{\ET}{\end{theo}}
\newcommand{\BK}{\begin{korr}}
	\newcommand{\EK}{\end{korr}}
\newcommand{\EQD}{\stackrel{d}{=}}
\newcommand{\BEX}{\begin{example}}
	\newcommand{\EEX}{\end{example}}

\newcommand{\BD}{\begin{de}}
	\newcommand{\ED}{\end{de}}

\newcommand{\BIT}{\begin{itemize}}
	\newcommand{\EIT}{\end{itemize}}
\newcommand{\BDI}{\begin{description}}
	\newcommand{\EDI}{\end{description}}

\newcommand{\BRM}{\begin{remark}}
	\newcommand{\ERM}{\end{remark}}

\newcommand{\BEL}{\begin{lem}}
	\newcommand{\EEL}{\end{lem}}

\newcommand{\nelem}[1]{{Lemma \ref{#1}}}

\newcommand{\netheo}[1]{{Theorem \ref{#1}}}

\newcommand{\equaldis}{\stackrel{{d}}{=}}

\begin{document}

\title{%Double
Finite-time ruin probability for correlated Brownian motions}

\author{Krzysztof D\c{e}bicki}
\address{Krzysztof D\c{e}bicki, Mathematical Institute, University of Wroc\l aw, pl. Grunwaldzki 2/4, 50-384 Wroc\l aw, Poland}
\email{Krzysztof.Debicki@math.uni.wroc.pl}

\author{Enkelejd  Hashorva}
\address{Enkelejd Hashorva, Department of Actuarial Science, %\\Faculty of Business and Economics\\
University of Lausanne,\\
UNIL-Dorigny, 1015 Lausanne, Switzerland
}
\email{Enkelejd.Hashorva@unil.ch}

\author{Konrad Krystecki}
\address{Konrad Krystecki, Department of Actuarial Science, %\\Faculty of Business and Economics\\
	University of Lausanne,\\
	UNIL-Dorigny, 1015 Lausanne, Switzerland and
	 Mathematical Institute, University of Wroc\l aw, pl. Grunwaldzki 2/4, 50-384 Wroc\l aw, Poland
}
\email{Konrad.Krystecki@unil.ch}

\bigskip

\date{\today}
 \maketitle

 {\bf Abstract:} Let $(W_1(s), W_2(t)), s,t\ge 0$ be a bivariate Brownian motion with standard Brownian motion marginals and constant correlation $\rho \in (-1,1)$ and define the joint survival probability of both supremum functionals $\pi_\rho(c_1,c_2; u, v)$ by
 %We derive exact asymptotics of
 $$\pi_\rho(c_1,c_2; u, v)=\pk*{\sup_{s \in [0,1]} \left(W_1(s)-c_1s\right)>u,\sup_{t \in [0,1]} \left(W_2(t)-c_2t\right)>v} ,$$
where $c_1,c_2 \inr$ and $u,v$ are given positive constants.
Approximation of $\pi_\rho(c_1,c_2; u, v) $ is of interest for the analysis of ruin probability in bivariate Brownian risk model as well as in the study of bivariate test statistics. In  this contribution we derive tight bounds for $\pi_\rho(c_1,c_2; u, v)$ in the case $\rho \in (0,1)$ and obtain precise approximations by letting $u\to \IF$ and taking $v= au$ for some fixed positive constant $a$ and $\rho \in (-1,1).$
 \COM{
 as $u \to \infty$, where $(W_1, W_2)$ is a correlated Brownian motion with correlation $\rho\in(-1,1)$,
$c_1,c_2\in \R$ and $a>0$.
The relation between model parameters $a$ and $\rho$ induces
several various cases,
both dimension-reduction and full-dimensional type of the asymptotics.
}



\section{Introduction }
Consider the Brownian risk model
$(R_1(s),R_2(t))$ of two insurance risk portfolios
\bqny{ R_1(s)= u+ c_1s - W_1(s), \quad R_2(t)= v+ c_2t - W_2(t),\quad s,t\ge0 ,
}
where the random process of accumulated claims $\left(W_1(s), W_2(t)\right), s,t\ge 0$
is assumed to be jointly Gaussian, the initial capitals are $u,v$
and the corresponding premium rates are $c_1,c_2$.
 In order to specify the model completely we need to give the joint law of $(W_1,W_2)$. In view of e.g.,  \cite{delsing2018asymptotics} (see also  \cite{SIM}) a natural  choice is
 to suppose that marginally $W_i's$ are standard Brownian motions with constant correlation $\rho\in(-1,1),$ i.e.
\BQN \label{BB}
(W_1(s),W_2(t))=(B_1(s), \rho B_1(t)+ \sqrt{1- \rho^2} B_2(t)), \quad s,t\ge 0,
\EQN
where $B_1,B_2$ are two independent standard Brownian motions. %and $\rho \in (-1,1)$.

The ruin probability of a single portfolio in the time horizon $ [0,T], T>0$ is given by (see e.g., \cite{MandjesKrzys})
\BQN \label{nuk}
{ \pi}_T(c_i;u):=\pk*{ \inf_{t\in [0,T]} R_i(t) < 0}&=& \pk*{\sup_{t\in [0,T]} (W_i(t)- c_i t)> u}\notag \\
&=& \Phi\left(-\frac u{ \sqrt{T}} -{c_i\sqrt{T}}\right)+
e^{-2c_iu}\Phi\left(- \frac u{\sqrt{T}} +{c_i\sqrt{T}}\right)
\EQN
for  $i=1,2,$ any $u\geq 0$ and with $\Phi(x)=1-\Psi(x)=\pk*{B_1(1) \le x}.$
\newline
Define next the component-wise ruin probability on $[0,T]$ by
%in time horizon $T>0$
%for bivariate risk process $(R_1(s),R_2(t)), s,t\ge0$ by
\begin{eqnarray*}
\pi_{T,\rho}(c_1,c_2;u,v)
&=&
\pk*{ \inf_{s \in [0,T]} R_1(s) <0 , \inf_{t\in [0,T]} R_2(t)< 0}=
\pk*{ \sup_{s \in [0,T]} W_1^*(s) > u, \sup_{t\in [0,T]} W_2^*(t)> v},
\end{eqnarray*}
where $W_i^* (s)=W_i(s)-c_i s.$ By the self similarity of Brownian motion,
without loss of generality we shall suppose that $T=1$ and set $\pi_{\rho}(c_1,c_2;u,v):=\pi_{1,\rho}(c_1,c_2;u,v).$ Clearly, for the special case $\rho=0$ we have the explicit formula
$$ \pi_{0}(c_1,c_2;u,v)=  { \pi}_1(c_1;u){ \pi}_1(c_2;v)$$
for any $u,v$.

$\pi_{T,\rho}$
has been investigated in \cite{KZ16,HE1998, Met10, ShaoWang13, RS06}.
In particular,  when $\rho\not=0$, in \cite{HE1998}[Thm 2.2]
there was derived a formula for
$$
%p(c_1,c_2;u,v):=
\pk*{ \sup_{s \in [0,T]} W_1^*(s) \le u, \sup_{t\in [0,T]} W_2^*(t) \le v},$$
which is given in terms of infinite-series and Bessel functions.
Representations given there are complex and do not
allow to observe the behaviour of $\pi_{T,\rho}$.
Therefore, in this contribution  we focus on
the exact estimates and bounds which
give more tractable view of the behaviour of $\pi_{T,\rho}$.
Infinite-time horizon analog of $\pi_{T,\rho}$
is studied  in \cite{DJR19,DJR20}, where both logarithmic and exact asymptotics
for $\pi_{\infty,\rho}(c_1,c_2;u,u)$, as $u\to\infty,$ was derived. We note that
due to infiniteness of time-interval in the model considered in \cite{DJR19,DJR20},
both the details of the proofs and the type of the asymptotics are different than in this contribution.

In   \cite{SIM} the simultaneous ruin probability
$$ \overline\pi_{\rho}(c_1,c_2;u,au)= \pk*{ \exists s \in [0,1]:  W_1^*(s) > u,  W_2^*(s)> au}, \quad  a \le 1$$
has been studied. Note that taking $a\le 1$ is no restriction in view of the symmetry of the model. Therein an upper bound for $\overline\pi_{\rho}$ is derived in terms of $p_{u,au}:=\pk{ W_1^*(1)> u, W_2^*(1)> au}$. Dealing with $\pi_{\rho}(c_1,c_2;u,au)$ is more difficult (apart from the case $\rho=0$). It turns out that an accurate upper bound can also be derived for $\pi_\rho$ if $\rho \in (0,1).$
\BT\label{th.upper}
If $\rho \in (0,1)$, then for all $u, v \ge 0$ %\footnote{ EH: Can $v$ be negative and $u$ positive?}
\begin{eqnarray}
p_{u,v}   \le \pi_{\rho}(c_1,c_2;u,v)\le
A(c_1,c_2) p_{u,v} ,
\label{BuB}
\end{eqnarray}
where $1/A(x,y)=\Psi( \max(0,\frac{y-\rho x}{\sqrt{1-\rho^2}} )  )
\Psi(\max(0,x) )$.% and $p_{u,v}=  \pk*{ W_1^*(1)> u , W_2^*(T)> v }$.
\ET
\COM{\BRM \label{negative} Notice that if $v<0$, then $$ \pi_{\rho}(c_1,c_2;u,v)={ \pi}_1(c_1;u)$$
because $\pk*{W_2^*(0) > v}=1.$ Therefore we can restrict the model to $a \ge 0.$
\ERM}

The upper bound above is given in terms of $p_{u,v}$ and the constant $A(c_1,c_2)$, which does not depend on $u$ nor $v$. This suggests that asymptotically, as $u\to \IF$
\bqn{ \label{A}
	 \pi_{\rho}(c_1,c_2;u,au) \sim C p_{u,au},}
 where $C>0$ is some constant and $f \sim g$ means $\lim_{u \to \IF}\frac{f(u)}{g(u)}=1$.
 Such a behaviour is already observed for the probability of simultaneous ruin in \cite{SIM}.
As we shall show  in the next section,
which contains main results of this paper, this statement
does not apply for all $\rho \in (-1,0)$.
It appears that it is useful to divide the problem into several cases that are determined by the position
and size of the area that dominates the exact asymptotics and lead to their separate forms; see Theorems \ref{Th0}, \ref{Th1} in Section \ref{s.main}.
In Section \ref{s.proofs} we determine the behaviour of the joint variance
of our process and together with the corresponding Pickands lemma we prove claim of the main theorems. In Appendix we present some proofs to lemmas used in previous section.

 \COM{In this contribution we prove that for $\rho >0$ such an approximation is valid and we determine $C$ explicitly. The other cases are more difficult and lead to different asymptotic behaviour of $\pi_\rho.$}
% If $<\in general this is not the case when $\rho <0$.   On the other side $\pi_{\rho}(c_1,c_2;u,v) \ge p_{u,v}$ for any $u$ and $v$ and any $\rho \in [-1,1]$.

\section{Main Results}\label{s.main}
For the choice $v= au, a>0$ the bounds in \eqref{BuB} are  asymptotically equal (up to some constant) if $u\to \IF$. This motivates the approximation of $\pi_{\rho}(c_1,c_2;u,au)$ as $u\to \IF$ to be discussed in this section.
\COM{
In this section we derive exact asymptotics of
\begin{eqnarray}
\pi_{\rho}(c_1,c_2;u,au)=\pk*{ \sup_{s \in [0,1]} W_1^*(s) > u, \sup_{t\in [0,1]} W_2^*(t)> au}
\label{ruin}
\end{eqnarray}
as $u\to\infty$, for the bivariate risk model specified in the introduction allowing
}
Below  $c_1, c_2$ are given constants and without loss of generality we suppose that $a\in(0,1]$. Recall that $W_i^*(t)=W_i(t)-c_i t, t \ge 0.$ We divide the obtained results on two scenarios:
(i) case $1 > \rho \ge a> 0$, when
one coordinate asymptotically dominates the other, leading to the reduction of dimension phenomena, and
(ii) - the remaining case, where both coordinates contribute to the asymptotics.

\subsection{Dimension-reduction case}
Suppose that $1 > \rho \ge a> 0$. It appears that
in this case the asymptotics of $\pi_{\rho}(c_1,c_2;u,au)$  as $u\to\infty$ is dominated by
the extremal behaviour of $W_1^*$, while $W_2^*$ contributes to the asymptotics
only by a constant.

\BT \label{Th0}
%\begin{enumerate}%[label=(\roman*)]
%\item[(i)]
(i) If $\rho > a > 0, $ then
%\bqn{
$	\pi_{\rho}(c_1,c_2;u,au) \sim { \pi_1}(c_1;u).
$%\]%\pk*{\sup_{s \in [0,1]}W_1^*(s) > u} \sim 2 \pk{W_1^*(1)> u}.
%}
%\item[(ii)]

(ii) If $\rho = a\in (0,1),$
 then
 %\bqn{
$	\pi_{\rho}(c_1,c_2;u,au) \sim  \Phi\left( \frac{\rho c_1-c_2 }{\sqrt{1-\rho^2}}\right)
{ \pi_1}(c_1;u).
%\pk*{ \sup_{s \in [0,1]} W_1^*(s) > u} \sim 2 \Phi\left( \frac{\rho c_1-c_2 }{\sqrt{1-\rho^2}}\right) \pk{W_1^*(1)> u}.
$%\]
% }
%\end{enumerate}
\ET
One can check that (\ref{A}) is satisfied under the assumptions considered in this section.
We recall that
${ \pi_1}(c_1;u) \sim 2 \pk{W_1^*(1)> u}$, as $u\to\infty$.

\subsection{Full-dimensional case}
Consider now  scenario complementary to the dimension-reduction case, i.e.
$\rho \in (-1,1)$  and $a\in  (\max(0,\rho), 1]$.
It appears that this case requires much deeper analysis divided on several subcases
which need separate approach leading to five different forms of the asymptotics.

Before presenting the main result of this section we introduce some useful notation. Let $\varphi_{s,t}$ be the probability density function (pdf) of $(W_1(s),W_2(t))$ and let $\varphi_{t}:=\varphi_{1,t}.$ Next define $\Sigma_{X}$ to be covariance matrix of random vector $X$ and denote by
%In order to To further simplify the most commonly used in this paper notation, denote
\BQN\label{Sst}
\Sigma_{ s,t}:= \Sigma_{(W_1(s),W_2(t))}=\begin{pmatrix}
s&  \rho \min(s,t) \\
\rho \min(s,t) &  t
\end{pmatrix}
\EQN
%to be
 the covariance matrix of $(W_1(s),W_2(t)).$ Finally, define for $t \in (0,1]$
\BQN \label{constantM}
M_{c_1,c_2,t}&=&(0,c_2)\Sigma^{-1}_{1,t}(1,a)^\top-(c_1,c_2)\left(\frac{1-2 \rho^2t} { t-\rho^2 t^2}\Sigma^{-1}_{1,t}-\frac{1}{t-\rho^2 t^2}\left( {\begin{array}{cc}
1 & -\rho \\
-\rho & 0 \\
\end{array} } \right)\right)(1,a)^\top
\EQN
and let
$ A_a= \frac{1}{4a}(1-\sqrt{a^2+8}).$
\BT \label{Th1}  Let   $\rho \in (-1,1), a\in  (\max(0,\rho), 1]$ and set $t^*=\frac{a}{\rho(2a\rho-1)},
\lambda_1 =\frac{1-a\rho}{1-\rho^2}, \lambda_2=  \frac{a-\rho}{1-\rho^2}$.\\
(i) If $\rho> A_a $, then
\bqny{
	\pi_{\rho}(c_1,c_2;u,au) \sim  C_1 u^{-2} \varphi_{1}(u+c_1 ,au+c_2),
}
where
$$ C_1= \int_{\R^2} \pk*{\exists_{s,t \in [0,\infty)}:\begin{array}{ccc}W_1(s) - s> x \\ W_2(t) - at>y\end{array}} e^{\lambda_1 x + \lambda_2 y} dx dy \in (0,\IF).$$
(ii) If $a <1$ and $\rho= A_a $, then
\bqny{ \pi_{\rho}(c_1,c_2;u,au) \sim  C_2 u^{-1} \varphi_{1}(u+c_1 ,au+c_2),
}
where
$C_2= \frac{2a}{\lambda_1}\frac{\sqrt{2 \pi}}{\sqrt{\tau}} \Phi\left(\frac{M_{c_1,c_2,1}}{\sqrt{\tau}}\right)e^{\frac{M_{c_1,c_2,1}^2}{2\tau}}, \tau=\frac{\rho^4  - 2 a \rho^5 - 3 a^2 \rho^2 + 3 a^2 \rho^4 + a^2 }{(1 - \rho^2)^3}>0.$\\
%\vfill
(iii) If $a=1, \rho= A_a=-\frac{1}{2}$, then
\bqny{
	\pi_{\rho}(c_1,c_2;u,u) \sim C_3 u^{-1}\varphi_{1}(u+c_1 ,u+c_2),
}
where $ C_3= \frac{\sqrt{2 \pi}}{\sqrt{\tau}}
\left( \Phi\left(\frac{M_{c_1,c_2,1}}{\sqrt{\tau}}\right)e^{\frac{M_{c_1,c_2,1}^2}{2\tau}}
+ \Phi\left(\frac{M_{c_2,c_1,1}}{\sqrt{\tau}}\right)e^{\frac{M_{c_2,c_1,1}^2}{2\tau}}\right), \tau=\frac{4}{3}. $
\\
%\vfill
(iv) If $a <1$ and $\rho< A_a $, then
\bqny{\pi_{\rho}(c_1,c_2;u,au) \sim  C_4 u^{-1}  \varphi_{t^*}(u+c_1 ,au+c_2t^*),
}
where $ C_4= 2a \frac{\sqrt{2 \pi}}{\sqrt{\tau}}\frac{1}{1-2a\rho} e^{\frac{M_{c_1,c_2,t^*}^2}{2\tau}}, \tau= -  \frac{ \rho^3 (1 - 2 a \rho)^4}{2a (1 - a \rho)}>0.$
\\
%\vfill
(v) If $a=1$ and $\rho< A_a = -\frac{1}{2} $, then
\bqny{\pi_{\rho}(c_1,c_2;u,au) \sim   C_5^{(1)}u^{-1}\varphi_{t^*}(u+c_1 ,u+c_2t^*) + C_5^{(2)}u^{-1}\varphi_{t^*}(u+c_1t^* ,u+c_2) ,
}
where $ C_5^{(1)}= 2 \frac{\sqrt{2 \pi}}{\sqrt{\tau}}\frac{1}{1-2\rho}e^{\frac{M_{c_1,c_2,t^*}^2}{2\tau}}, C_5^{(2)}= 2 \frac{\sqrt{2 \pi}}{\sqrt{\tau}}\frac{1}{1-2\rho}e^{\frac{M_{c_2,c_1,t^*}^2}{2\tau}}, \tau= -  \frac{ \rho^3 (1 - 2 \rho)^4}{2(1 - \rho)}>0.$
\ET

%For the considered in this section scenario, the play between $\rho$ and $a$ leads to five qualitatively different
%asymptotics as given in Theorem \ref{Th1}.
In  case (i) we still have that \eqref{A} holds.
%with constant
%$$ C_{a,\rho}= \int_{\R^2} \pk*{\exists_{s,t \in [0,\infty)}:\begin{array}{ccc}W_1(s) - s> x \\ W_2(t) - at>y\end{array}} e^{\lambda_1 x + \lambda_2 y} dx dy\in (0,\IF),    $$
%where $\lambda_1 =\frac{1-a\rho}{1-\rho^2} > 0$ and $\lambda_2=  \frac{a-\rho}{1-\rho^2}>0$.
Hovewer the claim of \netheo{th.upper} (and (\ref{A})) is not true for $\rho \in (-1, A_a]$ and $a\in (0,1]$.
%with
% $$ A_a= \frac{1}{4a}(1-\sqrt{a^2+8}).$$
It relates to the fact that when
 $\rho<0$ is relatively big compared to $a$ (in terms of absolute value),
 then it is less likely that the ruin occurs simultaneously.
 Hence the region that determines the asymptotics is separated from point $(1,1)$
 and the ruin is truly non-simultaneous. In those cases we can observe that
 \bqn{ \label{B}
	 \pi_{\rho}(c_1,c_2;u,au) \sim C up_{u,au}, \quad u\to \IF}
for some constant $C>0.$
In case (i) we have a similar constant to what appears in \cite{SIM}, and similarly to \cite{SIM} we cannot
calculate its exact value. However, notice that for cases (ii)-(v), the constants can be given explicitly
for particular $a,\rho.$

\section{Proofs}\label{s.proofs}

\subsection{Proof of Theorem \ref{th.upper}}
Given two independent standard Brownian motions $B_1,B_2$ let
\[S_1:=\sup_{t\in[0,1]}(B_1(t)-c_1t), \quad  S_2:=\sup_{t\in[0,1]}\left(\sqrt{1-\rho^2}B_2(t)-(c_2-\rho c_1)t\right).\]
Additionally, let $g_{1}(\cdot), g_{2}(\cdot)$ be probability density functions of $S_1$
and $\mathcal{X}:=\sqrt{1-\rho^2}B_2(1)-(c_2-\rho c_1)$,
respectively. Since $W_2(t)=\rho B_1(t)+\sqrt{1-\rho^2}B_2(t),$ we have for $u,v \ge 0 $
\begin{eqnarray}
\pi_{\rho}(c_1,c_2;u,v)
\le
\pk*{ S_1> u, \rho S_1+S_2> v}.
\end{eqnarray}
For any $c,u \in \mathbb{R}$ we obtain
\[
\pk*{\sup_{t\in[0,1]}(B(t)-ct) > u}\le\frac{\pk*{B(1)>u+c}}{\Psi\left(\max\left(0,c\right)\right)}
\]
see e.g., \cite{SIM,KoW18}.
Hence setting $1/A(x,y) =\Psi( \max(0,\frac{y-\rho x}{\sqrt{1-\rho^2}} )  )
\Psi(\max(0,x) )$  we have
\begin{eqnarray*}
	\lefteqn{
		\pk*{ S_1> u, \rho S_1+S_2> v}
	}\\
	&=&
	\int_u^\infty
	\pk*{S_2> v-\rho x}g_{1}(x) dx
	%}
	\\
	&\le&
	\frac{1}{\Psi\left( \max\left(0,\frac{(c_2-\rho c_1)}{\sqrt{1-\rho^2}} \right)  \right)}
	\int_u^\infty
	\pk*{\mathcal{X}> v-\rho x}g_{1}(x) dx\\
	&=&
	\frac{1}{\Psi\left( \max\left(0,\frac{(c_2-\rho c_1)}{\sqrt{1-\rho^2}} \right)  \right)}
	\pk*{ S_1> u, \rho S_1+\mathcal{X}> v}\\
	&=&
	\frac{1}{\Psi\left( \max\left(0,\frac{(c_2-\rho c_1)}{\sqrt{1-\rho^2}} \right)  \right)}
	\int_{-\infty}^\infty
	\pk*{ S_1>u, S_1>\frac{v-x}{\rho}}g_{2}(x)dx\\
	&\le&
	A(c_1,c_2)\int_{-\infty}^\infty \pk*{ B_1(1)-c_1>u, \rho(B_1(1)-c_1)>v-x}g_{2}(x)dx \\
	&=&
	A(c_1,c_2)\pk*{ W_1(1)> u+ c_1 , W_2(1)> v+ c_2},
\end{eqnarray*}
hence the claim follows.
\QED

\COM{
	Next, by symmetry we get
	\begin{eqnarray*}
		\pk*{ S_1> u, \rho S_1+S_2> v} &\le& A(c_2,c_1)\pk*{ W_1(T)> u+ c_1 T, W_2(T)> v+ c_2T}.
	\end{eqnarray*}
	This completes the proof.
	\QED
}
\subsection{Proof of Theorem \ref{Th0}}
First note that for any $u>0$
\BQN \label{lower}
\pi_{\rho}(c_1,c_2;u,au) &\ge& \pk*{ \exists t\in [0,1]:  W_1^*(t) > u, W_2^*(t)> au}.
\EQN
\underline{Case $0<a<\rho < 1$}. Notice that for any $u>0$
$$  \pi_{\rho}(c_1,c_2;u,au) \le  \pk*{ \sup_{s \in [0,1]} W_1^*(s)  > u}. $$
In view of \cite{SIM}[Thm 2.1] applied to the lower bound (\ref{lower}), we get
as $u\to\infty$
$$  \pi_{\rho}(c_1,c_2;u,au) \ge 2\pk*{ W_1^*(1)  > u}(1+o(1))=  \pk*{ \sup_{s \in [0,1]} W_1^*(s)  > u}(1+o(1))
$$
and hence
$$  \pi_{\rho}(c_1,c_2;u,au) \sim  \pk*{ \sup_{s \in [0,1]} W_1^*(s)  > u} =\pi_{1}(c_1;u). $$
\underline{Case $0<a=\rho < 1$}.
The  asymptotics of the lower bound follows again from \cite{SIM}[Thm 2.1] applied to (\ref{lower}):
\begin{eqnarray}
\pi_{\rho}(c_1,c_2;u,\rho u)
\ge
\Phi\left( \frac{\rho c_1-c_2 }{\sqrt{1-\rho^2}}\right)
{ \pi_1}(c_1;u)(1+o(1)), \ \ u\to\infty.\label{lower.1}
\end{eqnarray}
Setting  $h_u:=1-\frac{1}{\sqrt{u}}$ we have the following upper bound
\BQNY
\pi_{\rho}(c_1,c_2;u,\rho u) &\le&
\pk*{\exists_{ s,t \in [h_u,1]}: W_1^*(s)> u, W_2^*(t)>\rho u}+\pk*{\exists_{s \in [0,h_u]}: W_1^*(s)> u}\\
&&+\pk*{\exists_{ s \in [h_u,1], t \in [0,h_u]}: W_1^*(s)> u,W_2^*(t)> \rho u}.
\EQNY
Since $\rho>0$, it follows from (\ref{nuk}) that, for some $C > 0$

$$\pk*{\exists_{ s \in [0,h_u]}: W_1^*(s)> u} \le Ce^{-u\sqrt{u}}{ \pi_1}(c_1;u)(1+o(1)).$$
Additionally, for any $u>0$ we have
\BQNY
\pk*{\exists_{ s \in [h_u,1],t \in [0,h_u]}: W_1(s)> u,W_2(t)> \rho u}\\
&&\hspace*{-8cm} \le \pk*{\exists_{ s \in [h_u,1], t \in [0,h_u] }: \frac{b_1(s,t)B_1(s)+b_2(s,t)(\rho B_1(t) + \sqrt{1-\rho^2}B_2(t))}{b_1(s,t)+\rho b_2(s,t)}> u},
\EQNY
where $\vk{b}(s,t):=  \Sigma^{-1}_{s,t} (1,\rho)^\top.$
Since for all $u$ large
\BQN\sup_{s \in [h_u,1], t \in [0,h_u]}Var\left(\frac{b_1(s,t)B_1(s)+b_2(s,t)(\rho B_1(t) + \sqrt{1-\rho^2}B_2(t))}{b_1(s,t)
+\rho b_2(s,t)}\right) \sim 1-\frac{\rho^2}{1-\rho^2}\frac{1}{u}(1+o(1)),
\EQN
 then using \cite{Pit96}[Thm 8.1] we obtain  for some $C, \bar C$ positive and sufficiently large $u$
 \BQNY
\pk*{\exists_{ s \in [h_u,1],t \in [0,h_u]}: W_1(s)> u,W_2(t)> \rho u}
&\le& Ce^{-
	\frac{u^2}{2}\frac{1}{1-\bar{C}\frac{1}{u}}}\\
&=& %o\bar{C}e^{-\frac{\rho^2}{2(1-\rho^2)}u}
%\pk*{W_1^*(1)> u})
o({ \pi_1}(c_1;u))
, \ u\to\infty.
\EQNY
%\eE{where $\bar C_1, \bar  C_2$ are both positive.}
Using the above, we have that, as  $u\to\infty$
\bqny{
     \pi_{\rho}(c_1,c_2;u,\rho u)
&\le&
\pk*{\exists_{ s,t \in [h_u,1]}: W_1^*(s)> u, W_2^*(t)>\rho u} (1+o(1))   \\
& \le &
\left(\pk*{ \exists_{ s,t\in [h_u,1]}:\ W_1^*(s)> u, W_2^*(t)>\rho u, \forall_{ s \in [h_u,1]}: u+\frac{1}{\sqrt{u}}>W_1^*(s)}\right. \\
    &&+ \left.\pk*{ \exists_{ s\in [h_u,1]}:\ W_1^*(s)> u+\frac{1}{\sqrt{u}}}\right)(1+o(1)).
}
Due to $\left(\ref{nuk}\right),$ for some $C>0$ and suffifiently large $u$, we have
\BQNY
%\lefteqn{
\frac{\pk*{ \exists_{ s\in [h_u,1]}:\ W_1^*(s)> u+\frac{1}{\sqrt{u}}}}
{{ \pi_1}(c_1;u)}
&\le&
\frac{\pk*{ \exists_{ s\in [0,1]}:\ W_1^*(s)> u+\frac{1}{\sqrt{u}}}}
{\pk*{  W_1^*(1)> u}}
%}
\\
&=& \frac{\Phi\left(-(u+\frac{1}{\sqrt{u}}) -c_1\right)+e^{-2c_1(u+\frac{1}{\sqrt{u}})}
           \Phi\left(- (u+\frac{1}{\sqrt{u}})+c_1\right)}
{\Phi\left(-u -c_1\right)}\\
&\le& C e^{-\frac{(u+\frac{1}{\sqrt{u}})^2}{2}+\frac{u^2}{2}}\\
&=& C e^{-\sqrt{u}-\frac{1}{2u}}.% \quad u \to \infty.
\EQNY
Moreover, since $\rho>0$, then for $\bar{c_2}=c_2-\rho c_1$ %as $u\to\infty$
\bqny{
%     \pi_{\rho}(c_1,c_2;u,au)
 %   & = &
\lefteqn{
\pk*{ \exists_{ s,t\in [h_u,1]}:\ W_1^*(s)> u, W_2^*(t)>\rho u, \forall_{ s \in [h_u,1]}: u+\frac{1}{\sqrt{u}}>W_1^*(s)}}\\
%\pk*{\exists_{ s,t\in [h_u,1]}:B_1(s)-c_1s>u,\rho B_1(t)+\sqrt{1-\rho^2}B_2(t)-c_2t>\rho u}}\\
    & \le & \pk*{\exists_{ s,t\in [h_u,1]}:B_1(s)-c_1s>u,\rho (u+\frac{1}{\sqrt{u}})+\sqrt{1-\rho^2}B_2(t)-\bar{c_2}t>\rho u}\\
    & \le & \pk*{\exists_{ s\in [0,1]}:B_1(s)-c_1s>u}\pk*{\exists_{ t\in [h_u,1]}: \sqrt{1-\rho^2}B_2(t)-\bar{c_2}t>-\rho\frac{1}{\sqrt{u}}}\\
    & = & %\pk*{\exists_{ s\in [h_u,1]}:B_1(s)-c_1s>u}\pk*{\exists_{ t\in [h_u,1]}:B_2(t)- \bar{c_2}t/\sqrt{1-\rho^2}>0 }\\
    %& \le &
    \pk*{\exists s\in [0,1]:W_1^*(s)>u}
\Phi\left( \frac{\rho c_1-c_2 }{\sqrt{1-\rho^2}}\right)(1+o(1)),\ \ u\to\infty.
  %  \pk*{B_2(1)>\frac{\bar{c_2}}{\sqrt{1-\rho^2}}}.
}
Thus
$\pi_{\rho}(c_1,c_2;u,\rho u)\le \Phi\left( \frac{\rho c_1-c_2 }{\sqrt{1-\rho^2}}\right){ \pi_1}(c_1;u)(1+o(1)),
$
which combined with the asymptotic lower bound (\ref{lower.1}) completes the proof.\QED
\subsection{Proof of Theorem \ref{Th1}}
First we state several technical lemmas that are used in the proof.
In order to make the structure of the  proof more transparent, all proofs of the lemmas are deferred to Appendix.
Suppose that $a \in (\max(0,\rho),1]$ and recall  $\Sigma_{s,t}$ defined in \eqref{Sst}.
%for $s,t\in [0,1], \rho \in (-1,1)$ that
%$$\Sigma_{ s,t}= \begin{pmatrix}
%s&  \rho \min(s,t) \\
%\rho \min(s,t) &  t
%\end{pmatrix} $$
% is the covariance matrix of $(W_1(s),W_2(t))$.
Denote below for $\vk a = (1,a)^\top$
 $$ q_{\vk a}(s,t):= \vk a ^\top \Sigma^{-1}_{s,t} \vk a = \frac{t- 2 a \rho \min (s,t) + a^2s}{st- (\rho \min (s,t))^2 },$$
 $$\vk{b}(s,t):=  \Sigma^{-1}_{s,t} \vk a = \frac{1}{st- (\rho \min (s,t))^2}(t- a\rho \min(s,t), as  -\rho \min(s,t))^\top$$
 and set
\bqn{
	 q_{\vk a}^*(s,t)= \min_{ \vk x \ge \vk a} q_{\vk x}(s,t), \quad  q_{\vk a}^*= \min_{s,t \in [0,1]}  q_{\vk a}^*(s,t).
	}
It is well-known that $q_{\vk a}^*(s,t)$ captures the asymptotics of $\pk{ W_1^*(s)> u, W_2^*(t) > au}$, i.e., for any $s,t>0$ we have the following logarythmic asymptotics
\bqn{ \label{EXA2}
 \limit{u} \frac{1}{u^2} 	\log \pk{ W_1^*(s)> u, W_2^*(t) > au} =  - \frac{q_{\vk a}^*(s,t)}{2}.
}	
Moreover, by \cite{DEBICKIKOSINSKI}, we have
\bqn{ \label{EXA}
 \limit{u} \frac{1}{u^2} 	\log \pk{\exists_{s,t \in [0,1]} W_1^*(s)> u, W_2^*(t) > au} =  - \frac{q_{\vk a}^*}{2}.
}	
Note that for $a > \max(0,\rho)$ we have $\vk b(s,t) > (0,0)^\top.$	
Below we present a lemma that solves the problem of optimizing $q_{\vk a}^*(s,t).$

\BEL \label{optimal} %Set $t^*= \frac{a}{\rho(2a\rho-1)}$ for $\rho \in (-1,1), a \le 1$ given.
If
  $a=1, \rho<-\frac{1}{2}$ then function
  $q_{\vk a}^*(s,t)$ attains its minimum on $[0,1]^2$ at
$(s^*,t^*):= (1,\frac{a}{\rho(2a\rho-1)})$ and $(\bar{s}^*,\bar{t}^*):= (\frac{a}{\rho(2a\rho-1)}, 1).$
  %$(s_1^*,t_1^*):= (1,t^*)$ and $(s_2^*,t_2^*):= (t^*, 1).$
 For any other $a \in (\max(0,\rho),1]$ function $q_{\vk a}^*(s,t)$ attains its unique minimum on $[0,1]^2$ at
 $$(s^*,t^*):=\begin{cases}
 (1,\frac{a}{\rho(2a\rho-1)}) , & { if} \  \frac{a}{\rho(2a\rho-1)} \in [0, 1] \\
 (1,1), & otherwise. \end{cases}
  $$
\EEL
%Since only simple derivative calculations are required to prove the above lemma, the proof is omitted.

Hereafter we use notation for optimizers of
$\min_{s,t \in [0,1]}  q_{\vk a}^*(s,t)$ as introduced in Lemma \ref{optimal}.
Due to the symmetry of the case (v) of Theorem \ref{Th1} (note that $a=1$ in this case), in the rest of lemmas
presented below we focus only on the analysis of local properties of $(W^*_1,W^*_2)$
in the neighbourhood of point $(s^*,t^*)$, where $s^*=1$.

Let in the following $k_{u}=  1- \frac{(k-1)\Delta}{u^2},l_{u}=  t^*- \frac{(l-1)\Delta}{u^2}, u>0, \Delta>0$ and set
$$E_{u,k}=[(k+1)_u, k_u], \quad E_{u,k,l}=E_{u,k} \times E_{u,l}, \quad E= [-\Delta,0] \times [-\Delta,0].$$
Define also $$\chi_{u,k,l}(s,t):=(\chi_{1,u,k}(s),\chi_{2,u,l}(t)):
=u\left(W_1(\frac{s}{u^2}+k_u) - W_1(k_u) - c_1\frac{s}{u^2}, W_2(\frac{t}{u^2}+l_u )- W_2(l_u) - c_2\frac{t}{u^2}\right)$$
and let
\begin{eqnarray*}
S_{k_u,l_u}:=-(c_1,c_2)(\Sigma^{-1}_{k_u,l_u}-\Sigma^{-1}_{1,t^*})(u+c_1,au+c_2 )^\top
+(\frac{c_1(k-1)\Delta}{u^2},\frac{c_2(l-1)\Delta}{u^2})\Sigma^{-1}_{1,t^*}(u+c_1 k_u,au+c_2 l_u)^\top.%\label{S_ku}
\end{eqnarray*}
\BEL \label{singB}
Let $\rho \in (-1,1), a \in (\max(0,\rho),1], l,k \le \frac{u\log(u)}{\Delta}$ and  $\Delta>0 $ be given constants.
Then, as $u\to \IF$
\BQNY
\pk*{\exists_{(s,t) \in E_{u,k,l}}: \begin{array}{ccc}
	W_1^*(s)>u \\
	W_2^*(t)> au
	\end{array}} &\sim& I(\Delta) u^{-2} \varphi_{t^*}(u+c_1 ,au+c_2 t^*)
e^{S_{k_u,l_u}}
%e^{-((c_1,c_2)(\Sigma^{-1}_{k_u,l_u}-\Sigma^{-1}_{1,t^*})(u+c_1,au+c_2 )^\top)}\\
%&&\times  e^{(\frac{c_1(k-1)\Delta}{u^2},\frac{c_2(l-1)\Delta}{u^2})\Sigma^{-1}_{1,t^*}(u+c_1 k_u,au+c_2 l_u)^\top}
e^{-\frac{1}{2}u^2 (q_a(k_u,l_u)-q_a(1,t^*))},
\EQNY
$I(\Delta)= 	\begin{cases}
\int_{\R^2} \pk*{\exists_{s,t \in [0,\Delta]}:W_1(s) - s> x, W_2(t) - at>y} e^{\lambda_1 x + \lambda_2 y} dx dy &l_u = k_u\\
\int_{\R}\pk*{\exists_{t \in [0,\Delta]}:W_2(t) -\frac{a-\rho}{t^*-\rho^2}t > y}  e^{\lambda_2 y}  dy \int_{\R} \pk*{\exists_{s \in [0,\Delta]}:W_1(s) -s > x}e^{\lambda_1 x} dx  &l_u > k_u \\
\int_{\R} \pk*{\exists_{s \in [0,\Delta]}:W_1(s) -\frac{1-a\rho}{1 -\rho^2t^*}s > x}e^{\lambda_1 x} dx \int_{\R}\pk*{\exists_{t \in [0,\Delta]}:W_2(t) -\frac{a}{t^*}t > y}  e^{\lambda_2 y}  dy &l_u < k_u
\end{cases}$
\\and $\lambda_1 = \begin{cases}
\frac{t^*-a\rho}{t^*-\rho^2}, &l_u > k_u \\
\frac{1-a\rho}{1-\rho^2t^*}, &l_u < k_u \\
\frac{1}{t^*}\frac{1-a\rho}{1-\rho^2} &l_u = k_u
\end{cases}$, $\lambda_2 = \begin{cases}
\frac{a-\rho}{t^*-\rho^2 }, &l_u > k_u \\
\frac{a-\rho t^*}{t^*(1 -\rho^2t^*)}, &l_u < k_u \\
\frac{1}{t^*}\frac{a-\rho}{1-\rho^2}, &l_u = k_u
\end{cases}.$
\newline
Additionally
\begin{equation}
\label{Int}
 \lim_{u \to \infty}\sup_{l,k = O(u\log u)}\int_{\R^2}\pk*{\exists_{(s,t) \in E}:
    \begin{array}{ccc}
	\chi_{u,k,l}(s,t)>(x,y)
	\end{array}
    \Bigg{|}
    \begin{array}{ccc}
	W_1^*(k_u )= u - \frac{x}{u}  \\
	W_2^*(l_u )= au - \frac{y}{u}
	\end{array} }  e^{\lambda_1 x + \lambda_2 y} dxdy\\
< \infty.
\end{equation}
\EEL
The proof of Lemma \ref{singB} is derived mainly  by utilising the same idea as in the proof of
the classical Pickands lemma, see e.g., Lemma D.1 in \cite{Pit96}
or more recent contributions \cite{Rolski,DJR19}.
We note that finiteness of (\ref{Int}) is important in the proof of Theorem \ref{Th1}, where in order to evaluate sum over many small intervals we apply similar technique to the one used in \cite{DIEKER2005}[Lemma 2]. \\
Lemmas \ref{1dFinite}, \ref{2dFinite} deal with $\lim_{\Delta \to \IF} I(\Delta)$ for  $l_u \neq k_u$
and  $l_u=k_u$ (see \nelem{singB}) respectively.
The proof of \nelem{1dFinite} follows straightforwardly from \eqref{nuk}, % and hence is omitted,
while the proof of \nelem{2dFinite}
is largely the same as the proof of finiteness of
two-dimensional Piterbarg-type constants given in \cite{SIM}, see also \cite{DJR19}.
We omit the detailed calculations.
%is more complex and hence is deferred to Appendix.
\BEL Let in the following $B$ be a standard Brownian motion. \\ \label{1dFinite}
i) For any $b, c>0$ such that  $2b>c$  we have
%and $c\not=0$ we have
$$\lim_{\Delta \to \IF} \int_{\R} \pk{ \sup_{ t\in [0,\Delta]} (B(t) - bt)> x} e^{ cx} dx =
%\frac{1}{c} \mathcal{P}_W^{2b/c-1} \in (0,\infty),
\frac{1}{2b-c}+\frac{1}{c}.
$$
%where $ \mathcal{P}_W^{h}= \E{ \sup_{s\in [0,\IF ) } e^{  \sqrt{2} W(s)-  s - hs  }}$ is the Piterbarg constant corresponding to $W$ and $h>0$.
ii) For any $b>0$
$$
\lim_{\Delta \to\infty}\frac{1}{\Delta}\int_{\R} \pk{ \sup_{ t\in [0,\Delta]} (B(t) - bt)> x} e^{ 2bx} dx
= b.
$$
\EEL
Define for any $a > \max(0,\rho), A=\begin{pmatrix}1 & 0 \\ \rho & \sqrt{1-\rho^2} \end{pmatrix} ,\Sigma= AA^\top$, $\vk q \inr^2,$ $\Delta\in (0,\IF),$ and $\vk a =(1,a)^\top$
\BQNY
I(\Delta, \vk q):=\int_{\R^2}
\pk*{\exists_{\vk s\in [0,\Delta]^2 }: A \vk [\vk B(\vk s ) - \vk q \cdot \vk s] > \vk x}
e^{  \vk a ^\top \Sigma^{-1} \vk x} d
\vk x,
\EQNY
where $\vk B(\vk s)$ is a two-dimensional standard Brownian motion with independent components and $\vk a \cdot \vk b$ denotes component-wise multiplication of vectors. Note that $I( \Delta, A^{-1} \vk a)$ is the constant $I(\Delta)$ that appears in the case $k_u=l_u$ in \nelem{singB}.

\BEL \label {2dFinite} For any $a \in (\max(0,\rho),1]$ we have $I(\Delta, \vk q) \in (0,\IF)$ and
%Moreover, for $\vk q=A^{-1} \vk a$ we have further
\BQN \limit{\Delta} I(\Delta, A^{-1} \vk a) \in (0,\IF).
\EQN
\EEL
\COM{
In the following lemma we investigate Taylor expansion of the variance of $(W_1,W_2)$ process near the optimizing points.
We omit the straightforward proof of this lemma.
%Since proof of the lemma is based on straightforward calculations, hence the proof is omitted.
\BEL \label{Taylor} For any $a \in (\max(0,\rho),1]$ the following Taylor expansion holds
\BQNY \label{eqK1}
q_{\vk a}(s,t) - q_{\vk a} (1,t^*)& = &
\tau_1 (1-s)+  \tau_2 (1-t)+  \tau_3 (1-s)^2+ \tau_4(t^*-t)^2+\tau_5 (1-s)(t^*-t)\\
&&+  O((t^*-t)^3+ (1-s)^3),\ \ {\rm as}\ (s,t) \to (1,t^*).
\EQNY
(i) If $t^*=1$ and $t<s$ then
$$ \tau_1= \frac{(1 - a \rho)^2}{(1 - \rho^2)^2}, \quad \tau_2= \frac{-\rho^2 + 2 a \rho^3 + a^2 - 2 a^2 \rho^2}{(1 - \rho^2)^2},  $$
$$\tau_3=  \frac{(1 - a \rho)^2}{(1 - \rho^2)^3}, \quad
\tau_4= \frac{\rho^4  - 2 a \rho^5 - 3 a^2 \rho^2 + 3 a^2 \rho^4 + a^2 }{(1 - \rho^2)^3}, \quad \tau_5= -\frac{2 \rho^2 (1-  a \rho)^2}{(1 - \rho^2)^3}.  $$
(ii) If $t^*= \frac{a}{\rho(2a\rho-1)}$ and $t<s$ then
$$ \tau_1= (1 - 2 a \rho)^2, \quad \tau_2=0, \quad \tau_3= \frac{(1-2 a \rho)^3}{1-a \rho}, \quad \tau_4= -  \frac{ \rho^3 (1 - 2 a \rho)^4}{2a (1 - a \rho)} , \quad
\tau_5= -\frac{2 \rho^2 (1 - 2 a \rho)^3}{1 - a \rho}. $$
(iii) If $t^*=1$ and $s < t$ then
$$ \tau_1=\frac{1 - 2 \rho^2 - a^2 \rho^2 + 2 a \rho^3}{(1-  \rho^2)^2}, \quad
\tau_2=  \frac{(a - \rho)^2}{(1 - \rho^2)^2}, $$
$$   \tau_3=\frac{1 - 3 \rho^2 + 3 \rho^4 +  a^2 \rho^4 - 2 a \rho^5}{(1 - \rho^2)^3}, \quad
\tau_4= \frac{ (a - \rho)^2}{(1-  \rho^2)^3}, \quad
\tau_5=   -\frac{2(a - \rho)^2 \rho^2}{(1-  \rho^2)^3}.
$$
\EEL}
 The next lemma focuses on the asymptotic behaviour of $S_{k_u,l_u}$ appearing in the exponent in \nelem{singB}.
%We recall definition
\BEL \label {Sigma}
If $a \in (\max(\rho,0),1]$, then for $k_u>l_u, \Delta>0, k,l\le \frac{u\log(u)}{\Delta}$ and  $M_{c_1,c_2,t^*}$ given in \eqref{constantM}
$$S_{k_u,l_u} = M_{c_1,c_2,t^*} \frac{\Delta(l-1)}{u}+O\left(\frac{(l-1)^2\Delta^2}{u^3}\right)+O\left(\frac{(k-1)\Delta}{u}\right).$$
\EEL

Next for $C_1>0, C_2, i,j \in \mathbb{R}$ let
$$Q_{i,j}=\sum_{l=[i]}^{[j]} \frac{\sqrt{C_1}\Delta}{u} e^{C_2\frac{(l-1)\Delta}{u}-\frac{C_1}{2}\frac{(l-1)^2\Delta^2}{u^2}}, $$
where $[n]$ denotes the integer part of $n.$.
The following asymptotic result is used several times in the proof of the main theorem.
\BEL \label{STI} For $C_1>0, C_2, i,j \in \mathbb{R}$ we have as $u \to \IF$
$$Q_{1,u\log(u)}\sim \sqrt{2\pi}\Phi\left(\frac{C_2}{\sqrt{C_1}}\right)e^{\frac{C_2^2}{2C_1}}, \quad Q_{-u\log(u),u\log(u)}\sim \sqrt{2\pi}e^{\frac{C_2^2}{2C_1}}.$$
\EEL
Since the proofs of the above two lemmas follow by straightforward calculations, they are omitted.

\prooftheo{Th1} In the first part of the proof we
show that in order to determine the exact asymptotics of $\pi_\rho$ one can restrict the parameter set $[0,1]^2$
 to the area of size $\frac{\log(u)}{u}$ around the optimising points that were found in \nelem{optimal}.
 Then the proof is split into five cases;  in each case the contributing interval
 responsible for the asymptotics is different and a bit different argument has to be used. In first four cases there is one clear optimal point of function $q_{a}(s,t)$ and the asymptotics focuses around this point. In the last case there are two optimal points of function $q_{a}(s,t)$ and hence we treat that case differently.
 %and we need to restrict the intervals further to obtain the exact asymptotics.
 Recall that $ \vk b(s,t)= \Sigma^{-1} _{s,t }  (1,a) , $  	where $\Sigma_{s,t}$ is the covariance matrix of $(W_1(s), W_2(t))$. Since $a \in (\rho ,1]$ implies that $\vk b(s,t)$ has both components positive, then for any $u>0$ and $H \subset [0,1]$ we obtain an upper bound
 \bqny{
 \pk*{ \exists s,t\in H:\ W_1(s)> u, W_2(t)> au  }
% 	  & \le& \pk*{ \exists s,t\in H: \ b_1(s,t) W_1(s) + b_2(s,t) W_2(t)> u (b_1(s,t)+ a b_2(s,t))  }\\
 	  & \le& \pk*{ \exists s,t\in H:  \frac{ b_1(s,t) W_1(s) + b_2(s,t) W_2(t)} {b_1(s,t)+ a b_2(s,t)}> u  }.
}
The random field $Z(s,t)=\frac{ b_1(s,t) W_1(s) + b_2(s,t) W_2(t)} {b_1(s,t)+ a b_2(s,t)}$ has
variance function equal to  $1/ q_{\vk a}(s,t)= 1/q^*_{\vk a}(s,t)$ with $\vk a= (1,a)^\top$.  Consequently, by \nelem{optimal} we have
%\footnote{\kkk{K: it seems that up to case (v) we implicitly consider only the scenario when
%there is 1 optimizer (only $(1,t^*)$ appears at the beginning of the proof).
%If I am right, it has to be mentioned at the beginning of the proof.}}
$$ \sigma^2:= \sup_{s,t\in [0,1]} Var(Z(s,t))= \sup_{s,t \in [0,1]} \frac{1}{q_{\vk a}(s,t)} = \frac 1 { q_{\vk a} (s^*,t^*)}> 0.$$
Up to the proof of case (v), which we analyze separately, we suppose that %consider next for simplicity only the case
$s^*=1$ and let $H_\varepsilon=([1- \varepsilon, 1] \times [t^*-\varepsilon, t^*+ \kappa\varepsilon],$ where $\kappa=0$ if $t^*=1$ and $\kappa=1$ otherwise.  Define
$$F_u= [1-  \frac{\log u}u , 1] \times [t^*-\frac{\log u}u , t^*+ \kappa\frac{\log u}u  ]. \quad $$
Let next
$$N_u:=\floor{\frac{u\log(u)}{\Delta}}, \quad E_{u,m}^1:=[(m+1)_u, m_u],\quad E_{u,j}^2:=[(j+1)_u, j_u],$$
where $m_{u}=  1- \frac{(m-1)\Delta}{u^2}, j_{u}=  t^*- \frac{(j-1)\Delta}{u^2}, \Delta>0$.

In the first step, which is common for  cases (i)-(iv), we observe that
\BQN \label{star} \pi_{\rho}(c_1,c_2; u, au)
&=& 	\pk*{ \exists (s,t)\in  F_u:\ W_1^*(s)> u, W_2^*(t)> au  }(1+o(1))\\
&\le&
\sum_{k=1}^{N_u}\sum_{l=1-\kappa (N_u+1)}^{N_u}\pk*{\exists_{s \in E_{u,k}^1, t \in E_{u,l}^2}: W_1^*(s)>u, W_2^*(t)> au  }
(1+o(1)), \label{stars}
\EQN
as $u\to \IF$, where (\ref{star}) follows from \cite{Pit96}[Thm 8.1] and is proven in detail in Appendix,
while (\ref{stars}) is due to Bonferroni inequality.

\underline{Case (i): $\rho>\frac{1}{4a}(1-\sqrt{a^2+8}).$}
According to \nelem{optimal}
there is exactly one minimizer of
$q_{\vk a}^*(s,t)$  on $[0,1]^2$:  $(s^*,t^*)=(1,1)$.
Our aim is to prove that
$$\lim_{\Delta \to \IF}\lim_{u \to \IF}\frac{\pi_{\rho}(c_1,c_2;u, au)}{ \pk*{ \sup_{s \in [1-\frac{\Delta}{u^2},1]} W_1^*(s) > u, \sup_{t\in [1-\frac{\Delta}{u^2},1]}  W_2^*(t)> au}}=1.$$
%\begin{figure}[H]
%    \centering
%    \includegraphics[width=200pt]{Contrib1.png}
%\end{figure}
For $\Delta>0$, by \nelem{singB}, we have
\begin{eqnarray}
\pi_{\rho}(c_1,c_2;u, au)
&\ge&
\pk*{\exists_{s,t \in [1-\frac{\Delta}{u^2},1]}: W_1^*(s)>u, W_2^*(t)> au  }\nonumber\\
&\sim&
I(\Delta) u^{-2} \varphi_{1}(u+c_1 ,au+c_2)\  {\rm as}\ u\to\infty,\nonumber
\end{eqnarray}
where
$$I(\Delta)=
\int_{\R^2} \pk*{\exists_{s,t \in [0,\Delta]}:\begin{array}{ccc}W_1(s) - s> x \\ W_2(t) - at>y\end{array}} e^{\lambda_1 x + \lambda_2 y} dx dy< \IF, \quad
\lambda_1 =\frac{1-a\rho}{1-\rho^2}>0, \lambda_2=  \frac{a-\rho}{1-\rho^2}>0.$$

Using Taylor expansion we get as $u \to \IF$
$$u^2(q_a(k_u,l_u)-q_a(1,1))=\tau_1(k-1)\Delta+\tau_2(l-1)\Delta+o(\frac{1}{u}),$$
where
$\tau_1= \frac{(1 - a \rho)^2}{(1 - \rho^2)^2}>0,
\tau_2= \frac{-\rho^2 + 2 a \rho^3 + a^2 - 2 a^2 \rho^2}{(1 - \rho^2)^2}>0.$
 Implementing the above into \nelem{singB}, together with \eqref{stars}, we get  as $u\to\infty$
\BQNY
\pi_{\rho}(c_1,c_2;u, au)
&\le&
\sum_{k=1}^{N_u}\sum_{l=1}^{N_u}\pk*{\exists_{s \in E_{u,k}^1, t \in E_{u,l}^2}: W_1^*(s)>u, W_2^*(t)> au }(1+o(1))\\
& \le & \sum_{k=1}^{N_u}\sum_{l=1}^{N_u}e^{-\tau_1(k-1)\Delta}e^{-\tau_2(l-1)\Delta}
%\pk*{\exists_{s,t \in (1-\frac{\Delta}{u^2},1)}: W_1^*(s)>u, W_2^*(t)> au  }(1+o(1))
I(\Delta) u^{-2} \varphi_{1}(u+c_1 ,au+c_2)(1+o(1))
\\
& \le & \frac{1}{1-e^{-\tau_1\Delta}}\frac{1}{1-e^{-\tau_2\Delta}}
I(\Delta) u^{-2} \varphi_{1}(u+c_1 ,au+c_2)
%\pk*{\exists_{s,t \in (1-\frac{\Delta}{u^2},1)}: W_1^*(s)>u, W_2^*(t)> au  }
(1+o(1)).
\EQNY
Hence, letting $\Delta\to\infty$, and using that
by \nelem{2dFinite}
$$C_1:=\lim_{\Delta\to\infty}I(\Delta)=
\int_{\R^2} \pk*{\exists_{s,t \in [0,\infty)}:\begin{array}{ccc}W_1(s) - s> x \\ W_2(t) - at>y\end{array}}
e^{\lambda_1 x + \lambda_2 y} dx dy< \IF$$
we get
$$\lim_{u \to \IF}\frac{\pi_{\rho}(c_1,c_2;u, au)}
{C_1 u^{-2} \varphi_{1}(u+c_1 ,au+c_2)}=1.$$
\newline
\underline{Case (ii): $\rho=\frac{1}{4a}(1-\sqrt{a^2+8}).$}
According to \nelem{optimal}
there is exactly one minimizer of
$q_{\vk a}^*(s,t)$  on $[0,1]^2$:  $(s^*,t^*)=(1,1)$.
% Our aim is to prove that
%$$\lim_{\Delta \to \IF}\lim_{u \to \IF}\frac{\pi_{\rho}(c_1,c_2;u,au)}{ \pk*{ \sup_{s \in [1-\frac{\Delta}{u^2},1]} W_1^*(s) > u, \sup_{t\in [1-\frac{\log(u)}{u},1]}  W_2^*(t)> au}}=1.$$
%\begin{figure}[H]
%    \centering
%    \includegraphics[width=200pt]{Contrib3.png}
%\end{figure}
Using \eqref{stars} we have for any $\Delta>0$, as $u\to\infty$
\BQNY
\pi_{\rho}(c_1,c_2;u, au) %& \le &  \pk*{ \exists (s,t)\in F_u :\ W_1(s)> u, W_2(t)> au  }(1+o(1))\\
& \le &\sum_{k=1}^{N_u}\sum_{l=1}^{N_u}\pk*{\exists_{s \in E_{u,k}^1, t \in E_{u,l}^2}: W_1^*(s)>u, W_2^*(t)> au  }(1+o(1))
\EQNY
and also, by Bonferroni inequality
\BQN \label{triple}
\lefteqn{\pi_{\rho}(c_1,c_2;u, au)}\nonumber \\
& \ge &  \pk*{ \exists (s,t)\in F_u :\ W_1(s)> u, W_2(t)> au  }\nonumber \\
& \ge &\sum_{k=1}^{N_u}\sum_{l=1}^{N_u}\pk*{\exists_{s \in E_{u,k}^1, t \in E_{u,l}^2}: W_1^*(s)>u, W_2^*(t)> au  }\nonumber \\
&& -  \sum_{k=1}^{N_u}\sum_{l=1}^{N_u}\sum_{m=l+1}^{N_u} \mathbb{P}\Bigg{\{}\exists_{s \in E_{u,k}^1, t_1 \in E_{u,l}^2,t_2 \in E_{u,m}^2}: W_1^*(s)>u, W_2^*(t_1)> au, W_2^*(t_2)> au  \Bigg{\}}.
\EQN
Using Taylor expansion we have that
$$u^2(q_a(k_u,l_u)-q_a(1,1))=\tau_1 (k-1)\Delta+\tau_4 \frac{(l-1)^2\Delta^2}{u^2}+o(\frac{k^2}{u^2})+o(\frac{l^3}{u^4}),$$
where $\tau_1=\frac{(1 - a \rho)^2}{(1 - \rho^2)^2}>0$ and $\tau_4=\frac{\rho^4  - 2 a \rho^5 - 3 a^2 \rho^2 + 3 a^2 \rho^4 +  a^2 }{(1 - \rho^2)^3}>0.$
By \nelem{singB} and \nelem{Sigma} as $u \to \infty$ we have
\BQNY
\lefteqn{\sum_{k=1}^{N_u}\sum_{l=1}^{N_u}\pk*{\exists_{s \in E_{u,k}^1, t \in E_{u,l}^2}: W_1^*(s)>u, W_2^*(t)> au  }}\\
& \sim & I_1(\Delta) I_2(\Delta) u^{-2}\varphi_{t^*}(u+c_1 ,au+c_2)\sum_{k=1}^{N_u}\sum_{l=1}^{N_u}e^{-\frac{\tau_1}{2}(k-1)\Delta}e^{M_{c_1,c_2,t^*}\frac{(l-1)\Delta}{u}-\frac{\tau_4}{2}\frac{(l-1)^2\Delta^2}{u^2}} \\
& = & I_1(\Delta) I_2(\Delta) \frac{1}{1-e^{-\frac{\Delta\tau_1}{2}}} u^{-2}\varphi_{t^*}(u+c_1 ,au+c_2) \sum_{l=1}^{N_u} e^{M_{c_1,c_2,t^*}\frac{(l-1)\Delta}{u}-\frac{\tau_4}{2}\frac{(l-1)^2\Delta^2}{u^2}}  \\
& = & I_1(\Delta) \frac{I_2(\Delta)}{\Delta} \frac{1}{1-e^{-\frac{\Delta\tau_1}{2}}} u^{-1}\frac{1}{\sqrt{\tau_4}}\varphi_{t^*}(u+c_1 ,au+c_2) \sum_{l=1}^{N_u} \frac{\sqrt{\tau_4}\Delta}{u} e^{M_{c_1,c_2,t^*}\frac{(l-1)\Delta}{u}-\frac{\tau_4}{2}\frac{(l-1)^2\Delta^2}{u^2}},
\EQNY
where
$$I_1(\Delta)=\int_{\R} \pk*{\exists_{s \in [0,\Delta]}:W_1(s) -\frac{1-a\rho}{1-\rho^2}s > x}e^{\frac{1-a\rho}{1-\rho^2} x} dx, \quad
I_2(\Delta)=\int_{\R} \pk{ \sup_{ t\in [0,\Delta]} (W_2(t) - at)> x} e^{ 2ax} dx.$$\\
In view of  \nelem{1dFinite} $ \limit{\Delta} I_1(\Delta) = 2\frac{1-\rho^2}{1-a \rho}$ and $ \limit{\Delta} \frac{I_2(\Delta)}{\Delta} = a.$ Therefore by \nelem{STI}
\BQNY
\lim_{\Delta \to \IF}\lim_{u \to \IF}\frac{\sum_{k=1}^{N_u}\sum_{l=1}^{N_u}\pk*{\exists_{s \in E_{u,k}^1, t \in E_{u,l}^2}: W_1^*(s)>u, W_2^*(t)> au }}{2a \frac{1-\rho^2}{1-a \rho}\frac{\sqrt{2 \pi}}{\sqrt{\tau_4}} \Phi\left(\frac{M_{c_1,c_2,t^*}}{\sqrt{\tau_4}}\right)e^{\frac{M_{c_1,c_2,t^*}^2}{2\tau_4}} u^{-1}\varphi_{t^*}(u+c_1 ,au+c_2)}
& = & 1.
\EQNY
In order to complete the proof for this case, (\ref{triple}) needs to be shown to be asymptotically negligible,
which follows
by standard calculations  as in e.g., \cite{Pit96}. We defer those tedious calculations to Appendix.

\underline{Case (iii): $\rho=-\frac{1}{2}, a=1.$}
According to \nelem{optimal}
there is exactly one minimizer of
$q_{\vk a}^*(s,t)$  on $[0,1]^2$:  $(s^*,t^*)=(1,1)$.
%According to \nelem{optimal} $t^*=1$.
%Our aim is to prove that
%\BQNY
%\lim_{\Delta \to \IF}\lim_{u \to \IF}\frac{\pi_{\rho}(c_1,c_2;u,u)}{P_{u,\Delta}}&=&1,
%\EQNY
%where
%$$
%P_{u,\Delta}:=\pk*{ \exists_{s \in [1-\frac{\Delta}{u^2},1], t\in [1-\frac{\log(u)}{u},1]}\begin{array}{ccc}W_1^*(s)> u \\ W_2^*(t)>u\end{array}}+\pk*{ \exists_{s \in [1-\frac{\log(u)}{u},1], t\in [1-\frac{\Delta}{u^2},1]} \begin{array}{ccc}W_1^*(s)> u \\ W_2^*(t)>u\end{array}}.
%$$
%\begin{figure}[H]
%    \centering
%    \includegraphics[width=200pt]{Contrib4.png}
%\end{figure}
With \eqref{stars}, as $u\to \infty$
\BQNY
\pi_{\rho}(c_1,c_2;u, u) &\le& \sum_{k=1}^{N_u}\sum_{l=1}^{N_u-k}\pk*{\exists_{s \in E_{u,k+l}^1, t \in E_{u,l}^2}: W_1^*(s)>u, W_2^*(t)> u  }(1+o(1)) \\
&&+ \sum_{l=1}^{N_u}\sum_{k=1}^{N_u-l}\pk*{\exists_{s \in E_{u,k}^1, t \in E_{u,k+l}^2}: W_1^*(s)>u, W_2^*(t)> u  }(1+o(1)) \\
&&+ \sum_{k=1}^{N_u}\pk*{\exists_{s \in E_{u,k}^1, t \in E_{u,k}^2}: W_1^*(s)>u, W_2^*(t)> u  }(1+o(1))\\
&:=& \left(\sum_{k=1}^{N_u} S_{1,k} + \sum_{l=1}^{N_u} S_{2,l} + S_3 \right)(1+o(1)).
\EQNY
On the other hand, by Bonferroni inequality,
\BQNY
\lefteqn{\pi_{\rho}(c_1,c_2;u, u)}\\
&\ge& \pk*{ \exists (s,t)\in F_u :\ W_1(s)> u, W_2(t)> u  }\\
& \ge & \sum_{k=1}^{N_u}\sum_{l=1}^{N_u-k}\pk*{\exists_{s \in E_{u,k+l}^1, t \in E_{u,l}^2}: W_1^*(s)>u, W_2^*(t)> u  } \\
&&+ \sum_{l=1}^{N_u}\sum_{k=1}^{N_u-l}\pk*{\exists_{s \in E_{u,k}^1, t \in E_{u,k+l}^2}: W_1^*(s)>u, W_2^*(t)> u  } \\
&&+ \sum_{k=1}^{N_u}\pk*{\exists_{s \in E_{u,k}^1, t \in E_{u,k}^2}: W_1^*(s)>u, W_2^*(t)> u  }\\
&&-  \sum_{k=1}^{N_u}\sum_{l=1}^{N_u}\sum_{m=l+1}^{N_u} \mathbb{P}\Bigg{\{}\exists_{s \in E_{u,k}^1, t_1 \in E_{u,l}^2,t_2 \in E_{u,m}^2}: W_1^*(s)>u, W_2^*(t_1)> u, W_2^*(t_2)> u  \Bigg{\}}\\
&&-  \sum_{k=1}^{N_u}\sum_{l=1}^{N_u}\sum_{m=k+1}^{N_u} \mathbb{P}\Bigg{\{}\exists_{s_1 \in E_{u,k}^1, s_2 \in E_{u,m}^2,t \in E_{u,l}^2}: W_1^*(s_1)>u, W_2^*(s_2)> u, W_2^*(t)> u  \Bigg{\}}\\
&&- \sum_{k=2}^{N_u}\sum_{l=2}^{N_u}\pk*{\exists_{s_1 \in E_{u,l}^1, s_2 \in E_{u,0}^1, t_1 \in E_{u,0}^2, t_2 \in E_{u,k}^2}:
\begin{array}{ccc}
	W_1^*(s_1)>u \\
    W_1^*(s_2)>u \\
	W_2^*(t_1) > u \\
    W_2^*(t_2) > u
\end{array}
}\\
& := & \sum_{k=1}^{N_u} S_{1,k} + \sum_{l=1}^{N_u} S_{2,l} + S_3 - R_1 - R_2 - R_3.
\EQNY
\COM{K: Maybe we can delete the following (at this level of generality it doesn't give much insight):\\

Hence using \nelem{singB} and \eqref{Pick} we obtain for natural $k\le \frac{u\log(u)}{\Delta}$, as $u\to \infty$
$$
S_{1,k}\sim \sum_{l=1}^{N_u-k} e^{-4(l-1)}\pk*{\exists_{s \in E_{u,0}^1, t \in E_{u,k}^2}: W_1^*(s)>u, W_2^*(t)> u  } \sim  \pk*{\exists_{s \in E_{u,0}^1, t \in  E_{u,k}^2}:
\begin{array}{ccc}
W_1^*(s)>u \\
W_2^*(t)> u
\end{array} }.
$$
Similarly,
$$
S_{2,l} \sim \sum_{k=1}^{N_u-l} e^{-4(k-1)}\pk*{\exists_{s \in E_{u,l}^1, t \in E_{u,0}^2}: W_1^*(s)>u, W_2^*(t)> u  } \sim  \pk*{\exists_{s \in E_{u,l}^1, t \in E_{u,0}^2}:
\begin{array}{ccc}
 W_1^*(s)>u \\
 W_2^*(t)> u
\end{array} }
$$
and $S_3 \sim \pk*{\exists_{s \in E_{u,0}^1, t \in E_{u,0}^2}: W_1^*(s)>u, W_2^*(t)> u  }.$
\\
}
Observe that since $a=1,\rho=-\frac{1}{2},$ then for $k \le \frac{u\log(u)}{\Delta}, C>0,\Delta>0$
\BQN \label{Pick}
q_1(1-\frac{(k+C)\Delta}{u^2}, 1-\frac{C\Delta}{u^2})-q_1(1,1) = \frac{12C\frac{\Delta}{u^2}+4k\frac{\Delta^2}{u^4}-8Ck\frac{\Delta^2}{u^4}-12C\frac{\Delta^2}{u^4}}{(1-k\frac{\Delta}{u^2}-C\frac{\Delta}{u^2})(3-3C\frac{\Delta}{u^2}+k\frac{\Delta}{u^2})} \sim  \frac{4C\Delta}{u^2}
\EQN
as $u \to \infty,$  following similar calculations as in case (ii) we obtain
$$\lim_{\Delta \to \IF}\lim_{u \to \IF} \frac{\sum_{k=1}^{N_u}S_{1,k}}{ \sqrt{\frac{2 \pi}{{\tau_4}}} \Phi\left(\frac{M_{c_1,c_2,t^*}}{\sqrt{\tau_4}}\right)e^{\frac{M_{c_1,c_2,t^*}^2}{2\tau_4}}u^{-1}\varphi_{1}(u+c_1 ,u+c_2)}=1,$$
$$\lim_{\Delta \to \IF}\lim_{u \to \IF} \frac{ \sum_{l=1}^{N_u}S_{2,l}}{ \sqrt{\frac{2 \pi}{{\tau_4}}} \Phi\left(\frac{M_{c_2,c_1,t^*}}{\sqrt{\tau_4}}\right)e^{\frac{M_{c_2,c_1,t^*}^2}{2\tau_4}}u^{-1}\varphi_{1}(u+c_1 ,u+c_2)}=1,$$
where $\tau_4=\frac{\rho^4  - 2 \rho^5 - 3 \rho^2 + 3 \rho^4 +  1 }{(1 - \rho^2)^3}>0.$
Using \nelem{singB} and \eqref{Pick} we have
 $$\lim_{\Delta \to \IF}\lim_{u \to \infty}\frac{S_3}{\sum_{k=1}^{N_u}S_{1,k}}=0. $$
Now observe that for any $l>1$
\BQNY
R_3
\le
\sum_{k=2}^{N_u}\pk*{\exists_{s \in E_{u,0}^1, t_1 \in E_{u,0}^2, t_2 \in E_{u,k}^2}:
\begin{array}{ccc}
	W_1^*(s)>u \\
	W_2^*(t_1) > u \\
    W_2^*(t_2) > u
\end{array}
}.
\EQNY
With this observation and calculations similar as in case (ii) we have that $R_1, R_2$ and $R_3$ are asymptotically negligible, which completes the proof of case (iii).

\underline{Case (iv): $\rho<\frac{1}{4a}(1-\sqrt{a^2+8}).$}
%According to \nelem{optimal}, $t^*=\frac{a}{\rho(2a\rho-1)}.$
According to \nelem{optimal}
there is exactly one minimizer of
$q_{\vk a}^*(s,t)$  on $[0,1]^2$:  $(s^*,t^*)=(1,\frac{a}{\rho(2a\rho-1)})$.
%Our aim is to prove that
%$$\lim_{\Delta \to \IF}\lim_{u \to \IF}\frac{\pi_{\rho}(c_1,c_2;u, au)}{\pk*{ \sup_{s \in [1-\frac{\Delta}{u^2},1]} W_1^*(s) > u, \sup_{t\in [t^*-\frac{\log(u)}{u},t^*+\frac{\log(u)}{u}]}  W_2^*(t)> au}}=1.$$
%\begin{figure}[H]
%    \centering
%    \includegraphics[width=200pt]{Contrib2.png}
%\end{figure}
Using \eqref{stars} we have for any $\Delta>0$
\BQNY
\pi_{\rho}(c_1,c_2;u, au) %&\le& \pk*{ \exists (s,t)\in F_u :\ W_1(s)> u, W_2(t)> au  }(1+o(1)) \\
& \le &\sum_{k=1}^{N_u}\sum_{l=-N_u}^{N_u}\pk*{\exists_{s \in E_{u,k}^1, t \in E_{u,l}^2}: W_1^*(s)>u, W_2^*(t)> au  }(1+o(1))
\EQNY
and, by Bonferroni inequality
\BQN \label{tripleIV}
\pi_{\rho}(c_1,c_2;u, au) &\ge& \pk*{ \exists (s,t)\in F_u :\ W_1(s)> u, W_2(t)> au  }\nonumber \\
& \ge &\sum_{k=1}^{N_u}\sum_{l=-N_u}^{N_u}\pk*{\exists_{s \in E_{u,k}^1, t \in E_{u,l}^2}: W_1^*(s)>u, W_2^*(t)> au  }\nonumber \\
&&- \sum_{k=1}^{N_u}\sum_{l=-N_u}^{N_u}\sum_{m=l+1}^{N_u} \mathbb{P}\Bigg{\{}\exists_{s \in E_{u,k}^1, t_1 \in E_{u,l}^2,t_2 \in E_{u,m}^2}:W_1^*(s)>u, W_2^*(t_1)> au, W_2^*(t_2)> au  \Bigg{\}}.
\EQN
The rest of the proof follows by calculations similar to whose given in case (ii), with additional use of
the asymptotic symmetry of the behaviour of the components in the above summands around point $(1,t^*)$.

\underline{Case (v): $a=1, \rho<-\frac{1}{2}.$}
%In this unique case there are two completely symmetrical optimizing points for the variance of $(X_1(s),X_2(t).$
% function $$t^*=\frac{a}{\rho(2a\rho-1)}$ and there is a symmetrical optimising point with $s^*=\frac{a}{\rho(2a\rho-1)}.$
According to \nelem{optimal}, there are two minimizers of
$q_{\vk a}^*(s,t)$  on $[0,1]^2$:  $(s^*,t^*)= (1,\frac{1}{\rho(2\rho-1)})$ and
$(\bar{s}^*,\bar{t}^*)= (\frac{1}{\rho(2\rho-1)},1).$
Since $t^*=\bar{s}^*=\frac{1}{\rho(2\rho-1)}$, then in the remaining part of the proof,
in order to simplify notation, we replace $\bar{s}^*$ by $t^*$.
Denote
$$F_{1,u}=[1-\frac{\Delta \log(u)}{u},1] \times [t^*-\frac{\Delta \log(u)}{u},t^*+\frac{\Delta \log(u)}{u}],F_{2,u}=[t^*-\frac{\Delta \log(u)}{u},t^*+\frac{\Delta \log(u)}{u}] \times [1-\frac{\Delta \log(u)}{u},1].$$

%Our aim is to prove that
%\BQNY
%\lim_{\Delta \to \IF}\lim_{u \to \IF}\frac{\pi_{\rho}(c_1,c_2;u,u)}{
%\pk*{ \exists_{(s,t) \in F_{1,u}}: W_1^*(s) > u,  W_2^*(t)> u}
%+\pk*{ \exists_{(s,t) \in F_{2,u}}: W_1^*(s) > u,  W_2^*(t)> u}}=1.
%\EQNY
%\begin{figure}[H]
%    \centering
%   \includegraphics[width=200pt]{Contrib5.png}
%\end{figure}
Using symmetry of the optimizing points and the same idea as in \eqref{star} we get that as $u \to \IF$
\BQNY
\pi_{\rho}(c_1,c_2;u, u) &\le& \pk*{ \exists_{(s,t) \in F_{1,u}}: W_1^*(s) > u, W_2^*(t)>u}(1+o(1))\\
&&+\pk*{ \exists_{(s,t) \in F_{2,u}}: W_1^*(s) > u, W_2^*(t)>u}(1+o(1)).
\EQNY
On the other hand as $u \to \IF$
\BQN
\pi_{\rho}(c_1,c_2;u,u) &\ge& \pk*{ \exists_{(s,t) \in F_{1,u}}: W_1^*(s) > u, W_2^*(t)>u}+\pk*{ \exists_{(s,t) \in F_{2,u}}: W_1^*(s) > u, W_2^*(t)>u}\nonumber\\
&&- \pk*{ \exists_{(s,t) \in F_{1,u},(s',t') \in F_{2,u}}: W_1^*(s) > u, W_2^*(t)>u,W_1^*(s') > u, W_2^*(t')>u}.\label{quatro}
\EQN
Following the same calculations as in case (ii) we have that
\BQNY
\pk*{ \exists_{(s,t) \in F_{1,u}}: W_1^*(s) > u, W_2^*(t)>u}
&\sim& 2 \frac{\sqrt{2 \pi}}{\sqrt{\tau_4}} u^{-1} e^{\frac{M_{c_1,c_2,t^*}^2}{2\tau_4}}\frac{1-\rho^2 t^*}{1-\rho} \varphi_{t^*}(u+c_1 ,u+c_2t^*).
\EQNY
By the symmetry of this scenario with respect to sets $F_{1,u}$ and $F_{2,u}$, following case (ii), we get that
\BQNY
\pk*{ \exists_{(s,t) \in F_{2,u}}: W_1^*(s) > u, W_2^*(t)>u} %&=&
%\pk*{ \exists_{(s,t) \in F_{1,u}} W_1^*(s) > u, W_2^*(t)>u}\\
%\pk*{ \sup_{s \in [\frac{1}{\rho(2\rho-1)}-\frac{\Delta \log(u)}{u},\frac{1}{\rho(2\rho-1)}+\frac{\Delta \log(u)}{u}]} W_1^*(s) > %u, \sup_{t\in [1-\frac{\Delta}{u^2},1]}  W_2^*(t)> u}\\
&\sim& 2 \frac{\sqrt{2 \pi}}{\sqrt{\tau_4}} u^{-1} e^{\frac{M_{c_2,c_1,t^*}^2}{2\tau_4}}\frac{1-\rho^2 t^*}{1-\rho} \varphi_{t^*}(u+c_1t^* ,u+c_2).
\EQNY
Similarly to cases (ii)-(iv), (\ref{quatro}) needs to be shown to be asymptotically negligible. Since calculations are again standard and follow the ideas implemented already in cases (ii) and (iv), we omit those calculations.

\QED
\begin{alphasection}
\section{Appendix}

%%%%%%%%%%%%%%%%%%%%%%%%%%%%%%%%%%%%%%%%%%%%%%%%%%%%%%%%%%%%%%%%%%%%

\subsection{Proof of \eqref{star}}
We follow notation introduced in the proof of Theorem \ref{Th1}.

Since
$\sup_{s,t \in [0,1]^2 \setminus H_\varepsilon} Var(Z(s,t)) < \sigma^2,$ then for any $\varepsilon >0$ small enough  Borell-TIS inequality (see e.g.,  \cite{Pit96}) yields
$$ \pk{ \exists (s,t) \in [0,1]^2 \setminus H_\varepsilon :
	W_1(s)> u, W_2(t)> au} \le e^{-  r\frac{u^2 }{2  \sigma^2}    }   $$
for all sufficiently large $u$ and some $r>1$. Hence also
\bqny{ \pk{ \exists (s,t) \in [0,1]^2 \setminus H_\varepsilon :
			W_1^*(s)> u, W_2^*(t)> au}&\le &  e^{-  r^*\frac{u^2 }{2  \sigma^2}    }
}
for all large $u$ and some $r^*>1$. %Now for any $s\not= t, (s,t) \in H_\varepsilon \setminus F_u$   with
Recall
$F_u= [1-  \frac{\log u}u , 1] \times [t^*-\frac{\log u}u , t^*+ \kappa\frac{\log u}u  ]$.
%where $\kappa=0$ if $t^*=1$ and $\kappa=1$ otherwise.
Using Taylor expansion, for all $(s,t) \in H_\varepsilon \setminus F_u$ we have
\bqn{
	\sigma^2 -  Var(Z(s,t)) \ge \tau  \left(\frac{\log u}u\right)^2  \quad
}
for some $\tau>0$ and all $u$ large. We have
\BQNY
\pk*{ \exists (s,t)\in H_\varepsilon \setminus F_u:\ W_1(s)> u, W_2(t)> au  } &\le& \pk*{ \exists (s,t)\in H_\varepsilon \setminus F_u, s\not=t :\ W_1(s)> u, W_2(t)> au  } \\
&&+\pk*{ \exists (s,t)\in H_\varepsilon \setminus F_u, s=t :\ W_1(s)> u, W_2(t)> au  }.
\EQNY
Consequently, since $Z$ is a H\"older continuous random field and we can choose $\varepsilon>0$ such that
$Var(Z(s,t)) >0 $ for all $(s,t) \in H_\varepsilon$, then applying \cite{Pit96}[Thm 8.1] for some $c_1,C_1, C_2>0$ we have that
\bqny{
	\lefteqn{	\pk*{ \exists (s,t)\in H_\varepsilon \setminus F_u, s\not=t :\ W_1(s)> u, W_2(t)> au  }}\\
	& \le& \pk*{ \exists (s,t)\in H_\varepsilon \setminus F_u, s\not =t :  \frac{ Z(s,t)}{ \sqrt{Var(Z(s,t))}}\sqrt{Var(Z(s,t))} > u   }\\
	&\le & \pk*{ \exists (s,t)\in H_\varepsilon:   \frac{ Z(s,t)}{ \sqrt{Var(Z(s,t))}}  > u/ \sqrt{  \sigma^2 - \tau  \left(\frac{\log u}u\right)^2}   }\\
	&\le & C_1 u^{c_1} e^{-   \frac{u^2}{ 2 \sigma^2 - 2 \tau  \left(\frac{\log u}u\right)^2} }\\
	&\le & C_1 e^{- C_2  (\log u)^2 } e^{-   \frac{u^2}{ 2 \sigma^2 } }
}
and thus for $C_1',C_2'>0$
\bqny{
	\lefteqn{\pk*{ \exists (s,t)\in H_\varepsilon \setminus F_u, s\not=t:\ W_1^*(s)> u, W_2^*(t)> au  }}\\
	& \le& 	\pk*{ \exists (s,t)\in H_\varepsilon \setminus F_u:\ 	W_1(s)> u - c_1, W_2(t)> au - c_2    } \\
	&\le & C_1' e^{- C_2  (\log u)^2 } e^{-   \frac{u^2}{ 2 \sigma^2 } }.	
}
We have that
\bqny{
	\pk*{ \exists (s,t) \in H_\varepsilon :\ W_1^*(s)> u, W_2^*(t)> au  } & \ge & \pk*{  W_1^*(1) > u, W_2^*(t^*) > au }.%\\
	%	&\sim & \frac{1}{ \lambda_1 \lambda_2} u^{-2} \varphi_{t^*}( u+ c_1, au + c_2 t^*), \quad u\to \IF.
}
If  $t^*=1$ by  \cite{SIM}[Proposition 1.1]  for some $C^*$ positive
\bqny{
	\lefteqn{	\pk*{ \exists s=t,  0 \le  s \le  1- (\log u)/u :\ W_1^*(s)> u, W_2^*(t)> au  }}\\
	& \le & C^*\pk*{  W_1^*(1- (\log u)/u) > u, W_2^*(1- (\log u)/u ) > au }\\
	&=& o( \pk*{  W_1^*(1) > u, W_2^*(1) > au }), \quad u\to \IF .
}
Similarly, when $t^*< 1$ by \nelem{optimal} we have that $q^*(1,1)<q^*(1,t^*),$ hence
\bqny{
	\pk*{ \exists s=t , s \in [0,1] :\ W_1^*(s)> u, W_2^*(t)> au  }
	& \le & C^*\pk*{  W_1^*(1) > u, W_2^*(1) ) > au }\\
	&=& o( \pk*{  W_1^*(1) > u, W_2^*(t^*) > au }), \quad u\to \IF ,
}
where the last claim follows from \nelem{optimal} and \eqref{EXA}. Consequently, since
$$\pk*{ \exists (s,t)\in [0,1]^2 \setminus H_\varepsilon :\ W_1(s)> u, W_2(t)> au  }=o(\pk*{ \exists (s,t)\in H_\varepsilon \setminus F_u:\ W_1(s)> u, W_2(t)> au  } )$$
and
$$\pk*{ \exists (s,t)\in H_\varepsilon \setminus F_u:\ W_1(s)> u, W_2(t)> au  }=o(\pk*{ \exists (s,t)\in F_u:\ W_1(s)> u, W_2(t)> au  } ),$$
then \eqref{star} holds.
\QED

%%%%%%%%%%%%%%%%%%%%%%%

\subsection{Proof of negligibility of (\ref{triple}).}
Let $\phi$ denotes the density of $N(0,1)$ random variable. For any $0 \le l \le N_u$ using independence of increments of respective Brownian motions we have
\BQNY
\lefteqn{\sum_{k=1}^{N_u-l}\pk*{\exists_{s \in E_{u,1}^1, t_1 \in E_{u,k+l}^2, t_2 \in E_{u,l}^2}:
		\begin{array}{ccc}
			W_1^*(s)>u \\
			W_2^*(t_1) > a u \\
			W_2^*(t_2) > au
		\end{array}
}}\\
&= & \sum_{k=1}^{N_u-l}\int_{\R}\phi(u+c_1-\frac{x}{u})\\
&&\times\pk*{\exists_{s \in E_{u,0}^1, t_1 \in E_{u,k+l}^2, t_2 \in E_{u,l}^2}:
	\begin{array}{ccc}
		W_1^*(s)>u \\
		W_2^*(t_1) > a u \\
		W_2^*(t_2) > au
	\end{array}
	\Bigg{|}W_1(1)=u+c_1-\frac{x}{u}}dx \\
& = & \sum_{k=1}^{N_u-l}\int_{\R}\phi(u+c_1-\frac{x}{u})\pk*{\exists_{s \in E_{u,0}^1}:W_1(s)-W_1(1)+c_1(1-s) >\frac{x}{u} }\\
&&\times\pk*{\exists_{t_1 \in E_{u,k+l}^2, t_2 \in E_{u,l}^2}:
	\begin{array}{ccc}
		W_2^*(t_1) > a u \\
		W_2^*(t_2) > au
	\end{array}
	\Bigg{|}W_1(1)=u+c_1-\frac{x}{u}}dx \\
& = & \sum_{k=1}^{N_u-l}\int_{\R}\phi(u+c_1-\frac{x}{u})\pk*{\exists_{s \in E_{u,0}^1}:W_1(s)-W_1(1)+c_1(1-s) >\frac{x}{u} }\\
&&\times\pk*{\exists_{t_1 \in E_{u,k+l}^2, t_2 \in E_{u,l}^2}:
	X_{x,u}(t_1,t_2)>0}dx,
\EQNY
where $X_{x,u}(t_1,t_2)=(X_{1,x,u}(t_1),X_{2,x,u}(t_2))$ is a bivariate Gaussian process with \[\E {X_{x,u}(t_1,t_2)}=-\Bigg{(}\begin{array}{ccc}
-c_2t_1+\rho t_1 (c_1-\frac{x}{u}) \\
-c_2t_2+\rho t_2 (c_1-\frac{x}{u})
\end{array}\Bigg{)}+\Bigg{(}\begin{array}{ccc}
-(a-\rho t_1)u \\
-(a-\rho t_2)u
\end{array}\Bigg{)}\] and
\[\Sigma_{X_{x,u}(t_1,t_2)}=\left(\begin{array}{ccc}
t_1-\rho^2 t_1^2 & t_1-\rho^2 t_1 t_2 \\
t_1-\rho^2 t_1 t_2 & t_2-\rho^2 t_2^2
\end{array}\right).\]
Notice that $X_{x,u}(t_1,t_2)$ is H\"older continuous. Denote
$$S_0=\pk*{\exists_{t_1 \in E_{u,l}^2, t_2 \in E_{u,l}^2}:X_{x,u}(t_1,t_2)>0}, \quad S_1=\sum_{k=2}^{N_u-l}\pk*{\exists_{t_1 \in E_{u,k+l}^2, t_2 \in E_{u,l}^2}:X_{x,u}(t_1,t_2)>0},$$
$$S_2=\pk*{\exists_{t_1 \in (1-\frac{(l+2)\Delta}{u^2},1-\frac{(l+1+\frac{1}{\sqrt{\Delta}})\Delta}{u^2}), t_2 \in E_{u,l}^2}:X_{x,u}(t_1,t_2)>0},$$
$$S_3=\pk*{\exists_{t_1 \in(1-\frac{(l+1+\frac{1}{\sqrt{\Delta}})\Delta}{u^2},1-\frac{(l+1)\Delta}{u^2}), t_2 \in E_{u,l}^2}:X_{x,u}(t_1,t_2)>0}.$$
Observe that for (\ref{triple}) to be negligible it is enough to show that as $u \to \infty$
$$\frac{S_1+S_2+S_3}{S_0} \to 0.$$
We have
\BQNY
\pk*{\exists_{t_1 \in E_{u,k+l}^2, t_2 \in E_{u,l}^2}:X_{x,u}(t_1,t_2)>0}
&\le& \pk*{\exists_{t_1 \in E_{u,k+l}^2, t_2 \in E_{u,l}^2}:
	X_{1,x,u}(t_1)+X_{2,x,u}(t_2)>0} \\
& \le & \pk*{\exists_{t_1 \in E_{u,k+l}^2, t_2 \in E_{u,l}^2}:
	\frac{X_{1,x,u}(t_1)+X_{2,x,u}(t_2)}{\sigma_{k,u}}>0},
\EQNY
where $\sigma_{k,u}^2=\max_{t_1 \in E_{u,k+l}^2, t_2 \in E_{u,l}^2} \eta_{u}^2(t_1,t_2)$
and $\eta_{u}^2(t_1,t_2):=Var(X_{1,x,u}(t_1)+X_{2,x,u}(t_2)).$ Since on a set $t_1 \in E_{u,k+l}^2, t_2 \in E_{u,l}^2$ we have
$$\lim_{u \to \IF}t_1=\lim_{u \to \IF}t_2=1,$$
then as $u \to \IF$ we have on the same set $t_1 \in E_{u,k+l}^2, t_2 \in E_{u,l}^2$ that
$$\frac{\partial \eta_{u}^2(t_1,t_2)}{\partial t_1}=3-(2\rho^2+2t_2\rho^2) \sim 3-4\rho^2, \quad  \frac{\partial \eta_{k,u}^2(t_1,t_2)}{\partial t_2}=1-(2\rho^2+2t_1\rho^2) \sim 1-4\rho^2.$$
Since $\rho^2<\frac{1}{4}$ then above derivatives are positive for all large u. Hence $\eta_{k,u}^2(t_1,t_2)$ attains its maximum at $t_1^*=1-\frac{(l+k)\Delta}{u^2}, t_2^*=1-\frac{l\Delta}{u^2}.$ Consequently
$$\sigma_{k,u}^2=4-4\rho^2-\frac{1}{u^2}(4l\Delta+3k\Delta-8l\Delta \rho^2- 4k\Delta \rho^2)+O(\frac{1}{u^4}).$$
Denote $\mu_u:=\E{X_1(t_1^*)+X_2(t_2^*)}=2au+c_2t_1^*+c_2t_2^*-\rho (t_1^*+t_2^*)(u+c_1-\frac{x}{u})$. Using \cite{Pit96}[Thm 8.1], there exist constants $C, C_2>0$ such that
\BQN \label{s1}
S_1 &\le& \sum_{k=2}^{N_u-l} \pk*{\exists_{t_1 \in E_{u,k+l}^2, t_2 \in E_{u,l}^2}:\frac{X_{1,x,u}(t_1)+X_{2,x,u}(t_2)}{\sigma_{k,u}}>0} \nonumber\\
&\le& \sum_{k=2}^{N_u-l} C\frac{\mu_u}{\sigma_{k,u}} e^{-\frac{\mu_u^2}{2\sigma_{k,u}^2}} \nonumber\\
&=& \sum_{k=2}^{N_u-l} C \frac{\mu_u}{\sigma_{k,u}}e^{-\frac{\mu_u^2(4-4\rho^2+\frac{1}{u^2}(4l\Delta+3k\Delta-8l\Delta \rho^2- 4k\Delta \rho^2)+O(\frac{1}{u^4}))}{2((4-4\rho^2)^2+O(\frac{1}{u^4}))}} \nonumber\\
&\le&C \frac{\mu_u}{\sigma_{0,u}}e^{-\frac{\mu_u^2(4-4\rho^2+\frac{1}{u^2}(4l\Delta-8l\Delta \rho^2)+O(\frac{1}{u^4}))}{2((4-4\rho^2)^2+O(\frac{1}{u^4}))}}\sum_{k=2}^{N_u-l} e^{-C_2 k( \Delta+O(\frac{1}{u^2}))} \nonumber\\
&\le&C \frac{\mu_u}{\sigma_{0,u}}e^{-\frac{\mu_u^2(4-4\rho^2+\frac{1}{u^2}(4l\Delta-8l\Delta \rho^2)+O(\frac{1}{u^4}))}{2((4-4\rho^2)^2+O(\frac{1}{u^4}))}}\frac{e^{-C_2\Delta}}{e^{C_2\Delta}-1}.
\EQN
In the above we used the fact that $4l\Delta+3k\Delta-8l\Delta \rho^2- 4k\Delta \rho^2>0.$ Similarly we get that
\BQN \label{s2}
S_2 &\le&C \frac{\mu_u}{\sigma_{0,u}}e^{-\frac{\mu_u^2(4-4\rho^2+\frac{1}{u^2}(4l\Delta-8l\Delta \rho^2)+O(\frac{1}{u^4}))}{2(4-4\rho^2+O(\frac{1}{u^4}))}}e^{-C_2\sqrt{\Delta}}.
\EQN
Using \nelem{singB} and \nelem{STI} we have as $u \to \IF$
\BQN \label{s3}
\frac{S_3}{S_0} &\le & \frac{\pk*{\exists_{t_1 \in(1-\frac{(l+1+\frac{1}{\sqrt{\Delta}})\Delta}{u^2},1-\frac{(l+1)\Delta}{u^2})}:X_1(t_1)>0}}{\pk*{\exists_{t_1 \in(1-\frac{(l+2)\Delta}{u^2},1-\frac{(l+1)\Delta}{u^2})}:X_1(t_1)>0}} \nonumber\\
& = & \frac{\int_{\R} \pk*{\exists_{s \in [0,\sqrt{\Delta}]}:W_1(s) -\frac{1-a\rho}{1-\rho^2}s > x}e^{\frac{1-a\rho}{1-\rho^2} x} dx\int_{\R} \pk{ \exists_{ t\in [0,\sqrt{\Delta}]}: W_2(t) - at> x} e^{ 2ax} dx}{\int_{\R} \pk*{\exists_{s \in [0,\Delta]}:W_1(s) -\frac{1-a\rho}{1-\rho^2}s > x}e^{\frac{1-a\rho}{1-\rho^2} x} dx\int_{\R} \pk{ \exists_{ t\in [0,\Delta]}: W_2(t) - at> x} e^{ 2ax} dx} \nonumber\\
&=& \frac{\int_{\R} \pk{ \exists_{ t\in [0,\sqrt{\Delta}]}: W_2(t) - at> x} e^{ 2ax} dx}{\int_{\R} \pk{ \exists_{ t\in [0,\Delta]}:W_2(t) - at> x} e^{ 2ax} dx} \nonumber\\
&=& \frac{\sqrt{\Delta}}{\Delta}\frac{\int_{\R} \frac{1}{\sqrt{\Delta}}\pk{ \exists_{ t\in [0,\sqrt{\Delta}]}:W_2(t) - at> x} e^{ 2ax} dx}{\int_{\R} \frac{1}{\Delta} \pk{ \exists_{ t\in [0,\Delta]}:W_2(t) - at> x} e^{ 2ax} dx} \nonumber\\
&=&\frac{\sqrt{\Delta}}{\Delta}>0.
\EQN
Hence combination of (\ref{s1}), (\ref{s2}) and (\ref{s3}) leads to
$$\frac{S_1+S_2+S_3}{S_0} \le \frac{C}{\sqrt{\Delta}}+e^{-C k \sqrt{\Delta}}+\frac{e^{-C_2\Delta}}{e^{C_2\Delta}-1} \to 0, \ \Delta \to \infty$$
establishing the proof.
\QED
\end{alphasection}
\section{Supplementary Materials}
This section consists of supplementary technical proofs.
\subsection{Proof of \nelem{optimal}}

Since it will be needed in the proof of Eq. (3.4) below we consider $a \in [\rho ,1]$ for the derivation of \eqref{deri} below.
Supposing that  $s \le  t, t= cs, s,t\in [0,1],c\ge 1$ we have
 $$ q_{\vk a} (s,cs)=
 \frac{c- 2a  \rho  + a^2}{s(c- \rho ^2) }. $$
Since $c>a\rho$ and $a \ge \rho,$ then $c-2a\rho +a^2 = c-a\rho+a(a-\rho) >0$. Hence $q_{\vk a}(s,cs)$ is strictly decreasing in $s\le t$. Consequently,
 $$ \min_{0 \le s \le t \le 1} q_{\vk a}(s,t)= \min_{ z\in [0,1]} q_{\vk a}(z,1).$$
  Similarly, for $s \ge  t, s= ct, s,t\in [0,1],c\ge 1$ we have
  $$ q_{\vk a}(ct,t)=
  \frac{1- 2a  \rho  + a^2c}{t(c- \rho ^2) }. $$
Since $1>a\rho$ and $a \ge \rho$, we have
 $1-2a\rho +a^2c \ge 1- 2 a\rho+ a^2 =  1-a\rho+a(a-\rho)> 0$. Hence $q_{\vk a}(ct,t)$ is strictly decreasing in $t \le s$. Consequently, for any $g\in (0,1)$

 \BQN  \label{deri}
  \min_{s,t\in [0,1]\times [0,g]}  q_{\vk a}(s,t)= \min_{(z_1,z_2)\in \{1\}\times[0,1]\cup [0,1]\times\{g\}} q_{\vk a}(z_1,z_2).
 \EQN
Next we suppose that  $a>\rho$. By the definition of $\vk b(s,t)$ we have  that it has both components positive for any $s,t\in (0,1]$
and therefore $ q_{\vk a}^*(s,t) =  q_{\vk a}(s,t).$
For any $s,t$ positive such that $\vk b(s,t) $ has positive components we have that $\vk a^*= (1,a)^\top$, which follows from the solution of quadratic programming problem in \cite{Rolski}[Remark 5.1]. Hence
$$ q_{\vk a}^*=\min_{s,t \in [0,1]}  q_{\vk a}^*(s,t)=\min_{s,t \in [0,1]}  q_{\vk a}(s,t)=
\min \Bigl( \min_{z\in [0,1]}q_{\vk a}(z,1),
\min_{z\in [0,1]}q_{\vk a}(1,z ) \Bigr).$$
Calculating the derivatives we obtain for $z \in [0,1]$
$$\frac{d}{dz}q_{\vk a}(1,z)=\frac{z^2(\rho^2-2a\rho^3)+z(2a^2\rho^2)-a^2}{(z-\rho^2z^2)^2} $$
and
$$\frac{d}{dz}q_{\vk a}(z,1)=\frac{z^2(a^2\rho^2-2a\rho^3)+z(2\rho^2)-1}{(z-\rho^2z^2)^2}. $$
By setting the above derivatives to 0, it follows that all potential minimisation points of $q_{\vk a}^*(s,t)$
%from the boundary of $ [0,1] \times [0,1]$ are:
for $(s,t)\in [0,1] \times [0,1]$ are:
\begin{enumerate}
	\item $(s,t)=(1,1)$,
	\item $(s,t)=(\frac{1}{\rho(2\rho-a)},1)$,
	\item $(s,t)=(\frac{1}{a\rho},1)$,
	\item $(s,t)=(1,\frac{a}{\rho(2a\rho-1)})$,
	\item $(s,t)=(1,\frac{a}{\rho})$.
\end{enumerate}
It is easy to check, that for (1)-(4) $q_{\vk a}^*(s,t)>1.$ Since points in (3) and (5) do not belong to the boundary of $ [0,1] \times [0,1]$ for any values of $a,\rho$, then they can be excluded. Note that point in (2) belongs to the boundary only if point in (4) also belongs to the boundary and for those values of $a,\rho$ we have
$$q_{\vk a}(\frac{1}{\rho(2\rho-a)},1)=(2\rho-a)^2 >  (2a\rho-1)^2=q_{\vk a}(1,\frac{a}{\rho(2a\rho-1)}).$$
Hence point in (2) also cannot be an optimal point.
For $a=1, \rho<-\frac{1}{2}$ the symmetry shows as that the symmetrical points are the optimal points with the same minimal value. This completes the proof.
\QED

%%%%%%%%%%%%%%%%%%%%%%%%%%%%%%%%%%%%%%%%%%%%%%%%%%%%%

\section{Proof of \nelem{singB}}
Let $A_u:=\left\{\begin{array}{ccc}
W_1^*(k_u )= u - \frac{x}{u}  \\
W_2^*(l_u )= au - \frac{y}{u}
\end{array}\right\}.$ For all the cases we can write
\BQNY
\lefteqn{\pk*{\exists_{(s,t) \in E_{u,k,l}}: W_1^*(s)>u, W_2^*(t)> au  }}\\
&=&
\int_{\R^2}\pk*{\exists_{(s,t) \in E}:
	\begin{array}{ccc}
		W_1^*(\frac{s}{u^2}+ k_u )>u \\
		W_2^*(\frac{t}{u^2}+ l_u )> au
	\end{array}
	\Bigg{|}
	A_u} \\	
&&\times u^{-2}\varphi_{k_u,l_u}( u + c_1k_u - \frac{x}{u} ,  au + c_2l_u - \frac{y}{u}    )
dxdy\\
&=&
\int_{\R^2}\pk*{\exists_{(s,t) \in E}:
	\begin{array}{ccc}
		W_1^*(\frac{s}{u^2}+k_u ) - W_1(k_u ) +W_1(k_u)>u \\
		W_2^*(\frac{t}{u^2}+l_u )- W_2(l_u ) +W_2(l_u)> au
	\end{array}
	\Bigg{|}
	A_u } \\	
&&\times u^{-2}\varphi_{k_u,l_u}( u + c_1 k_u - \frac{x}{u} ,  au + c_2 l_u - \frac{y}{u}    )dxdy \\
&=&
\int_{\R^2}\pk*{\exists_{(s,t) \in E}:
	\begin{array}{ccc}
		\chi_{u,k,l}(s,t)>(x,y)
	\end{array}
	\Bigg{|}
	A_u } \\	
&&\times u^{-2}\varphi_{k_u,l_u}( u + c_1 k_u - \frac{x}{u} ,  au + c_2 l_u - \frac{y}{u}    )dxdy.
\EQNY
Furthermore, if $k_u \le l_u$, then
\BQNY
\lefteqn{\varphi_{k_u,l_u}( u + c_1 k_u - \frac{x}{u} ,  au + c_2 l_u - \frac{y}{u}    )}\\
&=&
\frac{1}{2\pi |\Sigma_{k_u,l_u}|}e^{-\frac{1}{2}(u+c_1 k_u-\frac{x}{u},au+c_2 l_u-\frac{y}{u})\Sigma^{-1}_{k_u,l_u}(u+c_1 k_u-\frac{x}{u},au+c_2 l_u-\frac{y}{u})^\top}\\
&\sim&
\frac{1}{2\pi |\Sigma_{1,t^*}|}e^{-\frac{1}{2}(u+c_1 k_u,au+c_2 l_u)\Sigma^{-1}_{k_u,l_u}(u+c_1 k_u,au+c_2 l_u)^\top} e^{\frac{l_u-a\rho k_u}{l_u k_u-\rho^2(k_u)^2} x + \frac{a-\rho}{l_u-\rho^2k_u} y} \\
\COM{&\sim&
\frac{1}{2\pi |\Sigma_{1,t^*}|}e^{-\frac{1}{2}(u+c_1 k_u,au+c_2 l_u)\Sigma^{-1}_{k_u,l_u}(u+c_1 k_u,au+c_2 l_u)^\top} e^{\lambda_1 x + \lambda_2 y} \\}
&\sim&
\frac{1}{2\pi |\Sigma_{1,t^*}|}e^{-\frac{1}{2} ((u+c_1,au+c_2t^*)(\Sigma^{-1}_{1,t^*}-\Sigma^{-1}_{1,t^*})(u+c_1,au+c_2t^* )^\top)} e^{\lambda_1 x + \lambda_2 y}\\
&&\times e^{-\frac{1}{2}(u+c_1,au+c_2t^*)\Sigma^{-1}_{k_u,l_u}(u+c_1,au+c_2t^* )^\top}e^{(\frac{c_1(k-1)\Delta}{u^2},\frac{c_2(l-1)\Delta}{u^2})\Sigma^{-1}_{k_u,l_u}(u+c_1 k_u,au+c_2 l_u)^\top} \\
&\sim& \varphi_{t^*}( u + c_1 ,  au + c_2 t^*) e^{(\frac{c_1(k-1)\Delta}{u^2},\frac{c_2(l-1)\Delta}{u^2})\Sigma^{-1}_{1,t^*}(u+c_1 k_u,au+c_2 l_u)^\top}\\
&&\times  e^{-\frac{1}{2} ((u+c_1,au+c_2t^*)(\Sigma^{-1}_{k_u,l_u}-\Sigma^{-1}_{1,t^*})(u+c_1,au+c_2t^* )^\top)}e^{\lambda_1 x + \lambda_2 y}\\
&=& \varphi_{t^*}( u + c_1 ,  au + c_2 t^*) e^{(\frac{c_1(k-1)\Delta}{u^2},\frac{c_2(l-1)\Delta}{u^2})\Sigma^{-1}_{1,t^*}(u+c_1 k_u,au+c_2 l_u)^\top}\\
&&\times  e^{-\frac{1}{2} u^2(q_{k_u,l_u}-q_{1,t^*})}e^{-((c_1,c_2t^*)(\Sigma^{-1}_{k_u,l_u}-\Sigma^{-1}_{1,t^*})(u+c_1,au+c_2t^* )^\top)}e^{\lambda_1 x + \lambda_2 y}.
\EQNY

Similarly, for $k_u>l_u$ we have
\BQNY
\varphi_{k_u,l_u}( u + c_1 k_u - \frac{x}{u} ,  au + c_2 l_u - \frac{y}{u}    ) &\sim& \varphi_{t^*}( u + c_1 ,  au + c_2 t^*) e^{(\frac{c_1(k-1)\Delta}{u^2},\frac{c_2(l-1)\Delta}{u^2})\Sigma^{-1}_{1,t^*}(u+c_1 k_u,au+c_2 l_u)^\top}\\
&&\times  e^{-\frac{1}{2} u^2(q_{k_u,l_u}-q_{1,t^*})}e^{-((c_1,c_2t^*)(\Sigma^{-1}_{k_u,l_u}-\Sigma^{-1}_{1,t^*})(u+c_1,au+c_2t^* )^\top)}e^{\lambda_1 x + \lambda_2 y}.
\EQNY
Next we investigate
$$I_u=\int_{\R^2}\pk*{\exists_{(s,t) \in E}:
	\chi^*_{u,k,l,x,y}(s,t) }  e^{\lambda_1 x + \lambda_2 y} dxdy, $$
where $\chi^*_{u,k,l,x,y}(s,t)=(\chi^*_{1,u,k,l,x,y}(s),\chi^*_{2,u,k,l,x,y}(t)):=\left(\chi_{u,k,l}(s,t)|A_u\right)$, $s,t \in [-\Delta,0].$ It appears that the play between $k_u$ and $l_u$ influences the next steps of argumentation, hence the rest of the proof is divided into three cases, $k_u<l_u, k_u=l_u,k_u>l_u.$

$(i)$
Suppose that $k_u=l_u$. Then
%\[
$\mathbb{E} \{\chi^*_{u,k,l,x,y}(s,t)\}=-\frac{1}{k_u}\Bigg{(}\begin{array}{ccc}
s(u+c_1k_u-\frac{x}{u})\\
t(au+c_2k_u-\frac{y}{u})
\end{array}\Bigg{)}
$%\]
and the covariance matrix of $\chi^*_{u,k,l,x,y}(s,t)$ is equal to
\BQNY
\Sigma_{\left(\chi^*_{u,k,l,x,y}(s,t)\right)}
&=&\Bigg{(}\begin{array}{ccc}
	s & \rho \min(s,t) \\
	\rho \min(s,t) & t
\end{array}\Bigg{)}-u^{-2}\Bigg{(}\begin{array}{ccc}
	s & \rho s \\
	\rho t & t
\end{array}\Bigg{)}
\Bigg{(}\begin{array}{ccc}
	k_u & \rho k_u \\
	\rho k_u & l_u
\end{array}\Bigg{)}^{-1}
\Bigg{(}\begin{array}{ccc}
	s & \rho t \\
	\rho s & t
\end{array}\Bigg{)}\\
&=&\Bigg{(}\begin{array}{ccc}
	s & \rho \min(s,t) \\
	\rho \min(s,t) & t
\end{array}\Bigg{)}-O\left(\frac{\log(u)}{u}\right)\Bigg{(}\begin{array}{ccc}
	s^2 & \rho^2 st \\
	(\rho \min(s,t))^2 & t^2
\end{array}\Bigg{)},s,t\in[0,\Delta].
\EQNY
Additionally since
$$
\chi^*_{u,k,l,x,y}(s_1,t_1)-\chi^*_{u,k,l,x,y}(s_2,t_2)\stackrel{d}= u\Big{(}W_1(\frac{s_1}{u^2}+k_u)-W_1(\frac{s_2}{u^2}+k_u),W_2(\frac{t_1}{u^2}+l_u)-W_2(\frac{t_2}{u^2}+l_u)|A_u\Big{)},
$$
then
%\BQNY
$\Sigma_{\left(\chi^*_{u,k,l,x,y}(s_1,t_1)-\chi^*_{u,k,l,x,y}(s_2,t_2)\right)}=\Sigma_{\chi^*_{u,k,l,x,y}(|s_1-s_2|,|t_1-t_2|)}.
$ %\EQNY
Using the above and the continuous mapping theorem we get, as $u \to \infty$
\BQNY
I_u & \sim & \int_{\R^2}\pk*{\exists_{s,t \in [0,\Delta]}:
	\begin{array}{ccc}
		W_1(s)-s>x \\
		W_2(t)-at>y
\end{array}}  e^{\lambda_1 x + \lambda_2 y} dxdy.
\EQNY
It order to justify the application of dominated convergence theorem we show finitness of (3.8). With $\lambda_1,\lambda_2>0$ we get for sufficiently large $u$
\BQN \label{fin}
I_u & = & \int_{\R^2}\pk*{\exists_{(s,t) \in E}:
	\begin{array}{ccc}
		\chi^*_{1,u,k,l,x,y}(s)>x\\
		\chi^*_{2,u,k,l,x,y}(t)>y
\end{array}}  e^{\lambda_1 x + \lambda_2 y} dxdy \nonumber\\
&\le & \int_{\R_-}\int_{\R_-}e^{\lambda_1 x + \lambda_2 y} dxdy +\int_{\R_-}\int_{\R_+}\pk*{\exists_{s \in [0,\Delta]}: \chi^*_{1,u,k,l,x,y}(s)>x}  e^{\lambda_1 x + \lambda_2 y} dxdy \nonumber \\
&&+\int_{\R_+}\int_{\R_-}\pk*{\exists_{t \in [0,\Delta]}: \chi^*_{2,u,k,l,x,y}(t)>y}  e^{\lambda_1 x + \lambda_2 y} dxdy \nonumber \\
&&+\int_{\R_+}\int_{\R_+}\pk*{\exists_{(s,t) \in [0,\Delta]}:\chi^*_{1,u,k,l,x,y}(s)+\chi^*_{2,u,k,l,x,y}(t)>x+y}  e^{\lambda_1 x + \lambda_2 y} dxdy \nonumber \\
&\le &\frac{1}{\lambda_1 \lambda_2}+\frac{1}{\lambda_2}\int_{\R_+}C_1e^{-C_2x^2} e^{\lambda_1 x} dx \\
&&+\frac{1}{\lambda_1}\int_{\R_+}C_1e^{-C_2y^2}  e^{\lambda_2 y} dy +\int_{\R_+}\int_{\R_+}C_1e^{-C_2(x+y)^2} e^{\lambda_1 x + \lambda_2 y} dxdy  < \infty, \nonumber
\EQN
where \eqref{fin} follows from \cite{Pit96}[Thm 8.1] with some constants $C_1,C_2>0$.
\newline
$(ii)$
Suppose that $k_u<l_u.$ Then the increments $W_1(s+k_u u^2) - W_1(k_u u^2),W_2(t+l_u u^2)- W_2(l_u u^2)$ are mutually independent. Hence
\BQNY
I_u & = &\int_{\R^2}\pk*{\exists_{s \in [0,\Delta]}: \chi^*_{1,u,x,y}(s)>x} \pk*{\exists_{t \in [0,\Delta]}: \chi^*_{2,u,x,y}(t)>y}  e^{\lambda_1 x + \lambda_2 y} dxdy,
\EQNY
where $\chi^*_{1,u,x,y}(s):=\chi_{1,u,k}(s)|A_u$ is a Gaussian process with
\[\mathbb{E} \{\chi^*_{1,u,x,y}(s)\}=\frac{1}{u(l_uk_u - \rho^2 k_u^2)}(sl_u - \rho^2 s k_u) (u+c_1-\frac{x}{u})- c_1\frac{s}{u},\]
\[Var\left(\chi^*_{1,u,x,y}(s)\right)=s-s^2\frac{l_u-\rho^2k_u}{u^2(k_ul_u-\rho^2k_u^2)}=s-\frac{s^2}{u^2k_u}\] and $\chi^*_{2,u,x,y}(t):=\chi_{2,u,l}(t)|A_u$ is a Gaussian process with
\[\mathbb{E}\{\chi^*_{2,u,x,y}(t)\}=\frac{1}{u(l_uk_u - \rho^2 k_u^2)}(k_u t (au+c_2-\frac{y}{u})-\rho k_u t (u+c_1-\frac{x}{u}))- c_2\frac{t}{u},\]
\[Var \left(\chi^*_{2,u,x,y}(t)\right)=t-t^2\frac{k_u}{u^2(k_ul_u-\rho^2k_u^2)}=t-\frac{t^2}{u^2(l_u-\rho^2k_u)}.\]
Moreover, for each $0\ge s>t\ge-\Delta$, $\chi^*_{1,u,x,y}(s)-\chi^*_{1,u,x,y}(t)$ is normally distributed with
\COM{\[
	\mathbb{E}\left(\chi^*_{1,u,x,y}(s)-\chi^*_{1,u,x,y}(t)\right)=
	\frac{1}{u(l_uk_u - \rho^2 k_u^2)}((s-t)l_u - \rho^2 (s-t) k_u) (u+c_1-\frac{x}{u})- c_1\frac{(s-t)}{u},
	\]}
\[Var \left(\chi^*_{1,u,x,y}(s)-\chi^*_{1,u,x,y}(t)\right)=(s-t)-\frac{(s-t)^2}{u^2k_u}
\]
while
$\chi^*_{2,u,x,y}(s)-\chi^*_{2,u,x,y}(t)$ is normally distributed with
\COM{\BQNY
	\lefteqn{\mathbb{E}\left(\chi^*_{2,u,x,y}(s)-\chi^*_{2,u,x,y}(t)\right)}\\
	&=&\frac{1}{u(l_uk_u - \rho^2 k_u^2)}(k_u (s-t) (au+c_2-\frac{y}{u})-\rho k_u (s-t) (u+c_1-\frac{x}{u}))- c_2\frac{(s-t)}{u},
	\EQNY}
\[Var \left(\chi^*_{2,u,x,y}(s)-\chi^*_{2,u,x,y}(t)\right)=(s-t)-\frac{(s-t)^2}{u^2(l_u-\rho^2k_u)}.
\]

Hence, using that $Var\left(\chi^*_{i,u,x,y}(s)-\chi^*_{i,u,x,y}(t)\right)\le 2|s-t|$ for all large enough $u$, we conclude that\\
$\chi^*_{1,u,x,y}(s)$ weakly converges to $W_1(s)-s$ and $\chi^*_{2,u,x,y}(t)$ weakly converges to $W_2(t)-\frac{a-\rho}{t^*-\rho^2}t$ in $C[0,\Delta].$

Next we prove the finitness of \eqref{Int} to justify the application of dominated convergence theorem. We have
\BQNY
I_u &\le& \int_{\R^2}\pk*{\exists_{s \in [0,\Delta]}: \chi^*_{1,u,x,y}(s)>x }\pk*{\exists_{t \in [0,\Delta]}: \chi^*_{2,u,x,y}(t)>y}  e^{\lambda_1 x + \lambda_2 y} dxdy\\
&\le& \int_{\R_-}\int_{\R_-}e^{\lambda_1 x + \lambda_2 y} dxdy \\
&&+ \int_{\R_-}\int_{\R_+}\pk*{\exists_{s \in [0,\Delta]}:
	\chi^*_{1,u,x,y}(s)>x}  e^{\lambda_1 x + \lambda_2 y} dxdy \\
&&+ \int_{\R_+}\int_{\R_-}\pk*{\exists_{t \in [0,\Delta]}:
	\chi^*_{2,u,x,y}(t)>y}  e^{\lambda_1 x + \lambda_2 y} dxdy \\
&&+\int_{\R_+}\int_{\R_+}\pk*{\exists_{s \in [0,\Delta]}:
	\chi^*_{1,u,x,y}(s)>x}\pk*{\exists_{t \in [0,\Delta]}:\chi^*_{2,u,x,y}(t)>y}  e^{\lambda_1 x + \lambda_2 y} dxdy.
\EQNY
Since $Var(\chi^*_{i,u,x,y}(s)-\chi^*_{i,u,x,y}(t))\le 2|s-t|$ for all large enough $u$ by \cite{Pit96}[Thm 8.1]
\BQNY
I_u &\le& \frac{1}{\lambda_1 \lambda_2}+\frac{1}{\lambda_2}\int_{\R_+}C_1e^{-C_2x^2} e^{\lambda_1 x} dx \\
&&+ \frac{1}{\lambda_1}\int_{\R_+}C_1e^{-C_2y^2}  e^{\lambda_2 y} dy +\int_{\R_+}\int_{\R_+}C_1e^{-C_2(x^2+y^2)} e^{\lambda_1 x + \lambda_2 y} dxdy  < \infty.
\EQNY
From the above it holds that \eqref{Int} is finite. Combining it with the weak convergence proven above and the dominated convergence theorem, we obtain that
\BQNY
\lim_{u \to \infty} I_u &=& \int_{\R^2}\pk*{\exists_{s \in [0,\Delta]}: W_1(s)-s>x }\pk*{\exists_{t \in [0,\Delta]}: W_2(t)-\frac{a-\rho}{t^*-\rho^2}t>y}  e^{\lambda_1 x + \lambda_2 y} dxdy.
\EQNY
\newline
$(iii)$ Suppose that $k_u>l_u$. Then the increments $W_1(s+k_u u^2) - W_1(k_u u^2), W_2(t+l_u u^2)- W_2(l_u u^2)$ are mutually independent. Hence we have
\BQNY
I_u & = &\int_{\R^2}\pk*{\exists_{s \in [0,\Delta]}: \xi^*_{1,u,x,y}(s)>x} \pk*{\exists_{t \in [0,\Delta]}: \xi^*_{2,u,x,y}(t)>y}  e^{\lambda_1 x + \lambda_2 y} dxdy,
\EQNY
where $\xi^*_{1,u,x,y}(s):=\chi_{1,u,k}(s)|A_u$ is a Gaussian process with
\[\mathbb{E} \{\xi^*_{1,u,x,y}(s)\}=\frac{1}{u(l_uk_u - \rho^2 l_u^2)}(s l_u (u+c_1-\frac{x}{u})-\rho s l_u (au+c_2-\frac{y}{u}))- c_1\frac{s}{u},\]
\[Var\left(\xi^*_{1,u,x,y}(s)\right)=s-s^2\frac{l_u}{u^2(k_ul_u-\rho^2l_u^2)}=s-\frac{s^2}{u^2(k_u-\rho^2l_u)}\] and $\xi^*_{2,u,x,y}(t):=\chi_{2,u,l}(t)|A_u$ is a Gaussian process with
\[\mathbb{E} \{\xi^*_{2,u,x,y}(t)\}=\frac{1}{u(l_uk_u - \rho^2 l_u^2)} t (k_u - \rho^2 l_u)(au+c_2-\frac{y}{u})- c_2\frac{t}{u},\]
\[Var \left(\xi^*_{2,u,x,y}(t)\right)=t-t^2\frac{k_u-\rho^2l_u}{u^2(k_ul_u-\rho^2l_u^2)}=t-\frac{t^2}{u^2l_u}.\]
Moreover, for each $0\ge s>t\ge-\Delta$, $\xi^*_{1,u,x,y}(s)-\xi^*_{1,u,x,y}(t)$ is normally distributed with
\[Var \left(\xi^*_{1,u,x,y}(s)-\xi^*_{1,u,x,y}(t)\right)=(s-t)-\frac{(s-t)^2}{u^2(k_u-\rho^2l_u)}
\]
and
$\xi^*_{2,u,x,y}(s)-\xi^*_{2,u,x,y}(t)$ is normally distributed with
\[Var \left(\xi^*_{2,u,x,y}(s)-\xi^*_{2,u,x,y}(t)\right)=|s-t|-\frac{(s-t)^2}{u^2l_u}.
\]

Hence, using that $Var \left(\xi^*_{i,u,x,y}(s)-\xi^*_{i,u,x,y}(t)\right)\le 2|s-t|$ for all u large enough,\\
$\xi^*_{1,u,x,y}(s)$ weakly converges to $W_1(s)-\frac{t^*(1-a\rho)}{t^*-\rho^2t^*}s$ and $\xi^*_{2,u,x,y}(t)$ weakly converges to $W_2(t)-\frac{a\rho(1-\rho t^*)}{t^*-\rho^2t^*}t$ in $C[0,\Delta].$
This leads to
\BQNY
\lim_{u \to \infty} I_u =\int_{\R^2}\pk*{\exists_{s \in [0,\Delta]}: W_1(s)-\frac{1-a\rho}{1-\rho^2}s>x }\pk*{\exists_{t \in [0,\Delta]}: W_2(t)-\frac{a}{t^*}t>y}  e^{\lambda_1 x + \lambda_2 y} dxdy.
\EQNY
The finitness of \eqref{Int} and the application of the dominated convergence theorem can be proven identically as in case (ii). This completes the proof.

\section{Proof of \nelem{1dFinite}}
Ad $i)$.  The proof follows straightforwardly from the fact that
$$\pk{ \sup_{ t\in [0,\infty)} (B(t) - bt)> x}=\min(1,e^{-2bx})$$ for $x \in \R$.\\
Ad $ii)$.
Note that, by self-similarity of Brownian motion and the change of variables
$y=2bx$, we have
\begin{eqnarray*}
	\int_{\R} \pk{ \sup_{ t\in [0,1]} (B(t) - bt)> x} e^{ 2bx} dx
	&=&
	\int_{\R} \pk{ \sup_{ t\in [0,2b^2T]} \left(B\left(\frac{t}{2b^2}\right) - \frac{t}{2b}\right)> x} e^{ 2bx} dx\\
	&=&
	\int_{\R} \pk{ \sup_{ t\in [0,2b^2T]} (\sqrt{2}B(t) - t)> 2bx} e^{ 2bx} dx\\
	&=&
	\frac{1}{2b}\int_{\R} \pk{ \sup_{ t\in [0,2b^2T]} (\sqrt{2}B(t) - t)> y} e^{ y} dy.
\end{eqnarray*}
Hence, using that
\[
\lim_{T\to\infty}\frac{1}{T}\int_{\R} \pk{ \sup_{ t\in [0,1]} (\sqrt{2}B(t) - t)> y} e^{ y} dy=1
\]
see e.g., \cite{Pit96}§ we complete the proof.
\QED

\section{Proof of \nelem{2dFinite}}
With $a \in (\max(0,\rho) ,1]$ define
$$ \vk \lambda = \Sigma^{-1} \vk a > \vk {0}, \quad \lambda_1= \frac{1- a \rho}{1- \rho^2}, \quad \lambda_2= \frac{a- \rho}{1- \rho^2}.$$
Since $\lambda_1,\lambda_2$ are positive, then   for any $\Delta>0$ %we have that
$$\int_{ x_1 \le 0, x_2 \le 0} v(\x)\, dx \in (0,\IF),$$
where
$$v(\x)=\pk*{\exists_{\vk s\in [0,\Delta]^2 }: A \vk [\vk B(\vk s ) - \vk q \cdot \vk s] > \vk x}
e^{  \vk a ^\top \Sigma^{-1} \vk x}.$$

If one of the coordinates of $\vk x$ is negative, then the integral reduces to one-dimensional case and it follows easily from Lemma 3.3 that this integral is also bounded for any $\Delta>0$.
The finiteness of $I(\Delta, \vk q) $ for $\Delta>0$ follows if we show further that
$$ \int_{x_1>0, x_2>0} v(\x)\, d \x $$ is finite, which can be established by using the fact that
$$\pk*{\exists_{\vk s\in [0,\Delta]^2 }: A \vk [\vk B(\vk s ) - \vk q \cdot \vk s] > \vk x}\le e^{- c \vk x^\top \vk x}$$
for some $c>0$. Note that for positive $\mu_1,\mu_2$ and $(X,Y)$ a bivariate random vector with finite moment generating function % %positive and $\mu_1,\mu_u_2 \inr$
\bqny{
	\int_{\R^2} \pk{ X> s, Y> t}   e^{\mu_1 s + \mu_2 t}ds dt &=& \frac1 {\mu_1 \mu_2}
	\E{ e^{ \mu_1 X+ \mu_2 Y}}.
}
Next, for $\Delta=n \in \mathbb{N}, \vk \lambda $ defined previously
\BQNY
\int_{\R^2} v(\x)  \, d\x &\le&
\sum_{i=0}^{n-1} \sum_{j=0}^{n-1}
\int_{\R^2}
\pk*{\exists_{\vk s\in [i,i+1]\times [j,j+1]} A \vk [\vk B(\vk s ) - (A^{-1} \vk a) \cdot \vk s] > \vk x}
e^{  \vk a ^\top \Sigma^{-1} \vk x} d
\vk x \\
&=& \frac{1}{\lambda_1 \lambda_2} \sum_{i=0}^{n-1} \sum_{j=0, j \neq i}^{n-1}\E{ e^{   \vk \lambda^\top \vk M(i,j)  }}+\frac{1}{\lambda_1 \lambda_2} \sum_{i=0}^{n-1} \E{ e^{   \vk \lambda^\top \vk M(i,i)  }},
\EQNY
where (below supremum is taken  component-wise)
\BQNY  \vk{M}(i,i) &=&\sup_{\vk s\in [i,i+1]\times [i,i+1]} A \vk [\vk B(\vk s ) - (A^{-1} \vk a) \cdot \vk s]\\
&=& \sup_{\vk s\in [i,i+1]\times [i,i+1]} A \vk [\vk B(\vk s )-  (B_1(i),B_2(i))^\top - (A^{-1} \vk a)   \cdot (\vk s- (i,i)^\top+ (i,i)^\top ) ] + A
(B_1(i),B_2(i))^\top \\
&\EQD   & \sup_{\vk s\in [0,1]\times [0,1]} A \vk [\vk B(\vk s ) - (A^{-1} \vk a)   \cdot \vk s ] -A(A^{-1} \vk a)\cdot (i,i)^\top + A
(B_1^*(i),B_2^*(i))^\top \\
&=& \vk Q- \vk a \cdot (i,i)^\top + A
(B_1^*(i),B_2^*(i))^\top
\EQNY
and for $j>i$ (case $i>j$ yields similar result)
\BQNY  \vk{M}(i,j) &=&\sup_{\vk s\in [i,i+1]\times [j,j+1]} A \vk [\vk B(\vk s ) - (A^{-1} \vk a) \cdot \vk s]\\
&=& \sup_{\vk s\in [i,i+1]\times [j,j+1]} A \vk [(B_1(s)-B_1(i),B_2(i+1)-B_2(i))^\top + (0,B_2(t)-B_2(j))^\top +\\
&&+ (0,B_2(j)-B_2(i+1))^\top - (A^{-1} \vk a) \cdot \vk (s-i,t-aj)^\top - (A^{-1} \vk a) \cdot \vk (i,aj)^\top+(B_1(i),B_2(i))^\top]\\
&\EQD& \sup_{\vk s\in [0,1]\times [0,1]} A \vk [\vk B(\vk s ) - (A^{-1} \vk a) \cdot \vk s]+\\
&&+ A[(B_1^*(i),(B_2^*(j)-B_2^*(i+1)+B_2^*(1)+B_2^*(i)))^\top]-(i,aj)^\top,
\EQNY
where $(B_1^*,B_2^*)$ is an independent copy of $(B_1,B_2)$ and $\equaldis$ stands for equality in law. Observe that
$$ A( B_1(i), B_2(j))^\top = (B_1(i), \rho B_1(i)+ \rho^* B_2(j))^\top. $$
Since $\lambda_1+ \lambda_2 \rho=1$ we have
$$ \vk \lambda^\top  A( B_1(i), B_2(j))^\top  =  (\lambda_1 +\lambda_2 \rho)B(i)
+ \lambda_2 \rho^* B_2(j) = B_1(i) + \lambda_2 \rho^* B_2(j)$$
implying  for some  $C$ positive
\BQNY   §\log\E{ e^{   \vk \lambda^\top \vk M(i,i)  }}
&=& \log \E{e^{ \vk \lambda^\top  \vk Q}}
+ \frac{i}{2}+ \frac{(a- \rho)^2}{2(1- \rho^2)} i -  \frac{1- a\rho}{1- \rho^2}i
- \frac{a-\rho}{1- \rho^2} ai \\
&=& \log \E{e^{ \vk \lambda^\top  \vk Q}}  - i \Bigl[ \frac{2- 2a \rho -(1- \rho^2)}{2(1-\rho^2)} \Bigr ]
- i\frac{a- \rho}{1- \rho^2}[2a-(a-\rho)] \\
&\le & \log \E{e^{ \vk \lambda^\top  \vk Q}}  - C i
\EQNY
and for $j>i$ for some $C_1>0, C_2$
\BQNY   \log \E{ e^{   \vk \lambda^\top \vk M(i,j)  }}
&=& \log \E{e^{ \vk \lambda^\top  \vk Q}}
+ \frac{i}{2}+ \frac{(a- \rho)^2}{2(1- \rho^2)} (j+2) -  \frac{1- a\rho}{1- \rho^2}i
- \frac{a-\rho}{1- \rho^2} aj \\
&=& \log \E{e^{ \vk \lambda^\top  \vk Q}}  - i \Bigl[ \frac{2- 2a \rho -(1- \rho^2)}{2(1-\rho^2)} \Bigr ]
- j\frac{a- \rho}{1- \rho^2}[2a-(a-\rho)] + C_2 \\
&\le & \log \E{e^{ \vk \lambda^\top  \vk Q}}  - C_1 (i+j)+ C_2,
\EQNY
hence the claim follows.
\QED

\section{Proof of \nelem{Sigma}}
For $k_u>l_u$ we have as $u \to \infty$
\BQNY
\Sigma^{-1}_{k_u,l_u}&=&\frac{1}{k_u l_u - \rho^2 l_u ^2} \left( {\begin{array}{cc}
		l_u & - \rho l_u \\
		- \rho l_u & k_u \\
\end{array} } \right)\\
&=& \frac{1}{(1-\frac{(k-1) \Delta}{u^2}) (t^*-\frac{(l-1) \Delta}{u^2}) - \rho^2 (t^*-\frac{(l-1) \Delta}{u^2}) ^2} \left( {\begin{array}{cc}
		(t^*-\frac{(l-1) \Delta}{u^2}) & - \rho (t^*-\frac{(l-1) \Delta}{u^2}) \\
		- \rho (t^*-\frac{(l-1) \Delta}{u^2}) & (1-\frac{(k-1) \Delta}{u^2}) \\
\end{array} } \right)\\
&=& \frac{1}{t^*-\rho^2 (t^*)^2 - \Delta \frac{(k-1)t^*+(l-1)-2\rho^2(l-1)t^*}{u^2}+ O(\frac{l^2+k^2}{u^4})} \left( {\begin{array}{cc}
		(t^*-\frac{(l-1) \Delta}{u^2}) & - \rho (t^*-\frac{(l-1) \Delta}{u^2}) \\
		- \rho (t^*-\frac{(l-1) \Delta}{u^2}) & (1-\frac{(k-1) \Delta}{u^2}) \\
\end{array} } \right)\\
&=& \left(\frac{1+  \frac{1}{u^2} \frac{\Delta [ (k-1)t^*+(l-1)-2 \rho^2(l-1)t^*]} { t^*-\rho^2 (t^*)^2}  + O(\frac{k^2+l^2}{u^4})}{t^*-\rho^2 (t^*)^2 }  \right) \left( {\begin{array}{cc}
		(t^*-\frac{(l-1) \Delta}{u^2}) & - \rho (t^*-\frac{(l-1) \Delta}{u^2}) \\
		- \rho (t^*-\frac{(l-1) \Delta}{u^2}) & (1-\frac{(k-1) \Delta}{u^2}) \\
\end{array} } \right)\\
&=& \Sigma^{-1}_{1,t^*}+\frac{D_{k,l}}{u^2}+O(\frac{k^2+l^2}{u^4})\left({\begin{array}{cc}
		1 & 1 \\
		1 & 1 \\
\end{array}}\right),
\EQNY
where
$$D_{k,l}=\frac{\Delta [ (k-1)t^*+(l-1)-2 \rho^2(l-1)t^*]} { t^*-\rho^2 (t^*)^2}\Sigma^{-1}_{1,t^*}-\frac{1}{t^*-\rho^2 (t^*)^2}\left( {\begin{array}{cc}
	\frac{(l-1) \Delta}{u^2} & -\frac{(l-1) \Delta}{u^2} \\
	-\frac{(l-1) \Delta}{u^2} & \frac{(k-1) \Delta}{u^2} \\
	\end{array} } \right).$$
The highest order term in terms of $k$ is $\frac{(k-1)\Delta}{u}.$ Further calculating the coefficient of $\frac{(l-1)\Delta}{u}$ by substituting $k=1$ we obtain
\BQNY
S_{1,l_u}&=&-(c_1,c_2)(\Sigma^{-1}_{1,l_u}-\Sigma^{-1}_{1,t^*})(u+c_1,au+c_2 )^\top+(0,\frac{c_2(l-1)\Delta}{u^2})\Sigma^{-1}_{1,t^*}(u+c_1,au+c_2 l_u)^\top\\
&=&-(c_1,c_2)\left(\frac{D_{1,l}}{u^2}+O\left(\frac{(l-1)^2\Delta^2}{u^4}\right)\left({\begin{array}{cc}
		1 & 1 \\
		1 & 1 \\
\end{array}}\right)\right)(u+c_1,au+c_2 )^\top\\
&&+(0,\frac{c_2(l-1)\Delta}{u^2})\Sigma^{-1}_{1,t^*}(u+c_1,au+c_2 l_u)^\top\\
&=&\frac{(l-1)\Delta}{u}(0,c_2)\Sigma^{-1}_{1,t^*}(1,a)^\top+O\left(\frac{(l-1)\Delta}{u^2}\right)-(c_1,c_2)\frac{D_{1,l}}{u^2}(u+c_1,au+c_2 )^\top+O\left(\frac{(l-1)^2\Delta^2}{u^3}\right)\\
&=&\frac{(l-1)\Delta}{u}(0,c_2)\Sigma^{-1}_{1,t^*}(1,a)^\top+O\left(\frac{(l-1)\Delta}{u^2}\right)\\
&&-\frac{(l-1)\Delta}{u}(c_1,c_2)\left(\frac{1-2 \rho^2t^*} { t^*-\rho^2 (t^*)^2}\Sigma^{-1}_{1,t^*}-\frac{1}{t^*-\rho^2 (t^*)^2}\left( {\begin{array}{cc}
		1 & -\rho \\
		-\rho & 0 \\
\end{array} } \right)\right)(1,a )^\top+O\left(\frac{(l-1)^2\Delta^2}{u^3}\right)\\
&=& M_{c_1,c_2,t^*} \frac{\Delta(l-1)}{u}+O\left(\frac{(l-1)^2\Delta^2}{u^3}\right).
\EQNY
\QED

%%%%%%%%%%%%%%%%%%%%%%%%%%%%%%%%%%%%%%%%%%%%%%%%%%%%%

\section{Proof of \nelem{STI}}
Notice that as $u \to \IF$
\BQNY
Q_{1,u\log(u)}&\sim&\int_{0}^{\IF}e^{-\frac{x^2}{2}+\frac{C_2}{\sqrt{C_1}}x}dx=\sqrt{2\pi}e^{\frac{C_2^2}{2C_1}}\int_{0}^{\IF}\frac{1}{\sqrt{2\pi}}e^{-\frac{x^2-2\frac{C_2}{\sqrt{C_1}}x+\frac{C_2^2}{C_1}}{2}}dx\\
&=&\sqrt{2\pi}e^{\frac{C_2^2}{2C_1}}\int_{0}^{\IF}\frac{1}{\sqrt{2\pi}}e^{-\frac{(x-\frac{C_2}{\sqrt{C_1}})^2}{2}}dx=\sqrt{2\pi}\Phi\left(\frac{C_2}{\sqrt{C_1}}\right)e^{\frac{C_2^2}{2C_1}}.
\EQNY
Similarly
\BQNY
Q_{-u\log(u),u\log(u)}&\sim&\int_{-\IF}^{\IF}e^{-\frac{x^2}{2}+\frac{C_2}{\sqrt{C_1}}x}dx=\sqrt{2\pi}e^{\frac{C_2^2}{2C_1}}\int_{-\IF}^{\IF}\frac{1}{\sqrt{2\pi}}e^{-\frac{x^2-2\frac{C_2}{\sqrt{C_1}}x+\frac{C_2^2}{C_1}}{2}}dx\\
&=&\sqrt{2\pi}e^{\frac{C_2^2}{2C_1}}\int_{-\IF}^{\IF}\frac{1}{\sqrt{2\pi}}e^{-\frac{(x-\frac{C_2}{\sqrt{C_1}})^2}{2}}dx=\sqrt{2\pi}e^{\frac{C_2^2}{2C_1}}.
\EQNY
\QED

\section{Proof of (3.4)}
In view of \eqref{deri}  we have with $g= h_u$

\COM{ For $s \le  t, t= cs, s,t\in [0,h_u],c\ge 1$ we have
 $$ q_{\vk a} (s,cs)=
 \frac{c- 2a  \rho  + a^2}{s(c- \rho ^2) }. $$
Since $c>a\rho$ and $a \ge \rho,$ then $c-2a\rho +a^2 = c-a\rho+a(a-\rho) >0$. Hence $q_{\vk a}(s,cs)$ is strictly decreasing in $s\le t$. Consequently,
 $$ \min_{0 \le s \le t \le 1} q_{\vk a}(s,t)= \min_{ z\in [0,h_u]} q_{\vk a}(z,1).$$
  Similarly, for $s \ge  t, s= ct, s,t\in [0,1],c\ge 1$ we have
  $$ q_{\vk a}(ct,t)=
  \frac{1- 2a  \rho  + a^2c}{t(c- \rho ^2) }. $$
Since $1>a\rho$ and $a \ge \rho$, we have
 $1-2a\rho +a^2c \ge 1- 2 a\rho+ a^2 =  1-a\rho+a(a-\rho)> 0$. Hence $q_{\vk a}(ct,t)$ is strictly decreasing in $t \le s$. Consequently
}
  $$ \min_{s\in [0,1],t \in [0,h_u]} q_{\vk a}(s,t)= \min_{(z_1,z_2)\in \{1\}\times[0,1]\cup [0,1]\times\{h_u\}} q_{\vk a}(z_1,z_2).$$
  Furthermore $$\frac{\partial}{\partial s}q_{\vk a}(s,h_u)<0, s \le 1, \quad \frac{\partial}{\partial t}q_{\vk a}(1,t)<0, t < h_u$$
  implying that
 $$ \min_{s\in [0,1],t \in [0,h_u]} q_{\vk a}(s,t)= q_{\vk a}(1,h_u).$$

   \COM{Notice that point $(1,\frac{a}{\rho})=(1,1)$ from the proof of Lemma 3.1 belongs to the boundary of $[0,1]^2$ in this case and is the optimal point with $\frac{\partial}{\partial t}q_{\vk a}(1,t)|_{t=1}=0.$} The random field $Z(s,t)=\frac{ b_1(s,t) W_1(s) + b_2(s,t) W_2(t)} {b_1(s,t)+ \rho b_2(s,t)},s,t\in [0,1]$ has
variance function equal to  $1/ q^*_{\vk a}(s,t)$ with $\vk a= (1,\rho)^\top$.  The first component of  $\vk b(s,t)$ is positive and the second component is $a s - \rho \min(s,t)= a(s - \min(s,t))$ which is equal to zero for all $s\le t$ and is positive for $s> t$. Hence the solution of the quadratic programming problem $q_{\vk a}^*(s,t)$ is $(1, \rho)=(1,a)$ and thus $q_{\vk a}^*(s,t)= q_{\vk a}(s,t)$. Consequently,
% > (0,0)^\top$ for $s\neq t$ and $\vk b(s,t)\ge (0,0)^\top$ for $s= t$ where $s,t\in [0,1]$. Hence for any $u>0$
$$\sup_{s\in [0,1],t \in [0,h_u]} Var(Z(s,t))= \sup_{s\in [0,1],t \in [0,h_u]} \frac{1}{q_{\vk a}^*(s,t)} =
\frac{1}{ \inf_{s\in [0,1],t \in [0,h_u]}  q_{\vk a}(s,t)} =\frac 1 { q_{\vk a} (1,h_u)}.$$
Next, for $h_u= 1 - 1/ \sqrt{u}$
\BQNY
1-q_{\vk a} (1,h_u) &=&1-\frac{h_u-2\rho^2(h_u)+\rho^2}{h_u-\rho^2h_u^2}
%&=&\frac{2\rho^2(1-c)-\rho^2(1-c)^2-\rho^2}{1-\rho^2 - (c(1-2\rho^2)+c^2\rho^2)}\frac{1-\rho^2 + c(1-2\rho^2)+c^2\rho^2}{1-\rho^2 + %c(1-2\rho^2)+c^2\rho^2}\\
%&=&\frac{c^2\rho^2(\rho^2-1)+c^3(2\rho^4-\rho^2)-c^4\rho^4}{(1-\rho^2)^2-c^2(1-2\rho^2)^2-2c^3\rho^2(1-2\rho^2)-c^4\rho^4}\\
\sim- \frac{1}{u} \frac{\rho^2}{1-\rho^2}(1+o(1))
\EQNY
establishing the proof.
%Hence by the asymptotics of $q_{\vk a} (1,h_u)$ as $u \to \IF$ we obtain
%$$\sup_{s\in [0,1],t \in [0,h_u]} Var(Z(s,t))=1-\frac{\rho^2}{1-\rho^2}\frac{1}{u}(1+o(1)).$$
\QED

%%%%%%%%%%%%%%%%%%%%%%%%%%%%%%%%%%%%%%%%%%%%%%%%%%%%%

\section{Proof of case (iv) of Theorem 2.2}
Recall that
\BQNY  \pi_{\rho}(c_1,c_2;u, au) &\le& \sum_{k=1}^{N_u}\sum_{l=-N_u}^{N_u}\pk*{\exists_{s \in E_{u,k}^1, t \in E_{u,l}^2}: W_1^*(s)>u, W_2^*(t)> au  } \\
& + & \pk*{ \exists (s,t)\in [0,1]^2 \setminus F_u :\ W_1(s)> u, W_2(t)> au  } \\
& = & \sum_{k=1}^{N_u}\sum_{l=-N_u}^{N_u}\pk*{\exists_{s \in E_{u,k}^1, t \in E_{u,l}^2}: W_1^*(s)>u, W_2^*(t)> au  } (1+o(1)).
\EQNY
For any $u>0$
\BQNY
\pi_{\rho}(c_1,c_2;u, au) &\ge& \pk*{ \exists (s,t)\in F_u :\ W_1(s)> u, W_2(t)> au  } \\
& \ge &\sum_{k=1}^{N_u}\sum_{l=-N_u}^{N_u}\pk*{\exists_{s \in E_{u,k}^1, t \in E_{u,l}^2}: W_1^*(s)>u, W_2^*(t)> au  } \\
& - & \sum_{k=1}^{N_u}\sum_{l=-N_u}^{N_u}\sum_{m=l+1}^{N_u} \mathbb{P}\Bigg{\{}\exists_{s \in E_{u,k}^1, t_1 \in E_{u,l}^2,t_2 \in E_{u,m}^2}:W_1^*(s)>u, W_2^*(t_1)> au, W_2^*(t_2)> au  \Bigg{\}}.
\EQNY
Using Taylor expansion we have
$$u^2(q_a(k_u,l_u)-q_a(1,1))=\tau_1(k-1)\Delta+\tau_4\frac{(l-1)^2\Delta^2}{u^2}+o(\frac{k^2}{u^2})+o(\frac{l^3}{u^4}),$$
where $\tau_1=(1 - 2 a \rho)^2>0, \tau_4=-\frac{ \rho^3 (1 - 2 a \rho)^4}{2a (1 - a \rho)}>0.$ Using Lemma 3.2, Lemma 3.5 and the symmetry of the sum, we get as $u \to \infty$
\BQNY
\lefteqn{\sum_{k=1}^{N_u}\sum_{l=-N_u}^{N_u}\pk*{\exists_{s \in E_{u,k}^1, t \in E_{u,l}^2}: W_1^*(s)>u, W_2^*(t)> au  }}\\
& \sim & Iu^{-2}\varphi_{t^*}(u+c_1 ,au+c_2t^*) \sum_{k=1}^{N_u}\sum_{l=-N_u}^{N_u}e^{-\tau_1(k-1)\Delta}e^{M_{c_1,c_2,t^*}\frac{l\Delta}{u}-\frac{\tau_4}{2}\frac{l^2\Delta^2}{u^2}},
\EQNY
where $$I=\int_{\R} \pk{ \sup_{ s\in [0,\Delta]} (W_1(s) - \frac{1-a\rho}{1-\rho^2 t^*}s)> x} e^{ \frac{1-a\rho}{1-\rho^2 t^*}x} dx  \int_{\R}\pk{ \sup_{ t\in [0,\Delta]} (W_2(t) - \frac{a}{t^*}t)> x} e^{ 2\frac{a}{t^*}x} dx. $$
Using Lemma 3.3 with Lemma 3.6, we get as $u \to \IF$
\BQNY
\lefteqn{\sum_{k=1}^{N_u}\sum_{l=-N_u}^{N_u}\pk*{\exists_{s \in E_{u,k}^1, t \in E_{u,l}^2}: W_1^*(s)>u, W_2^*(t)> au  }}\\
& \sim & 2t^*u^{-1}\frac{1}{\sqrt{\tau_4}}\frac{1}{(1-e^{-\tau_1\Delta})}\varphi_{t^*}(u+c_1 ,au+c_2t^*) \frac{1-\rho^2 t^*}{1-a\rho}\\
&& \times  \int_{\R}\frac{1}{\Delta}\pk{ \sup_{ t\in [0,\Delta]} (W_2(t) - \frac{a}{t^*}t)> x} e^{ 2\frac{a}{t^*}x} dx \sum_{l=-N_u}^{N_u}\frac{\sqrt{\tau_4}\Delta}{u}e^{M_{c_1,c_2,t^*}\frac{l\Delta}{u}-\frac{\tau_4}{2}\frac{l^2\Delta^2}{u^2}}\\
& \sim & 2t^*u^{-1}\frac{1}{\sqrt{\tau_4}} \frac{1-\rho^2 t^*}{1-a\rho} \frac{a}{t^*} \varphi_{t^*}(u+c_1 ,au+c_2t^*)\sum_{l=-N_u}^{N_u}\frac{\sqrt{\tau_4}\Delta}{u}e^{M_{c_1,c_2,t^*}\frac{l\Delta}{u}-\frac{\tau_4}{2}\frac{l^2\Delta^2}{u^2}}\\
& \sim &2a \frac{\sqrt{2 \pi}}{\sqrt{\tau_4}} u^{-1} e^{\frac{M_{c_1,c_2,t^*}^2}{2\tau_4}}\frac{1-\rho^2 t^*}{1-a\rho} \varphi_{t^*}(u+c_1 ,au+c_2t^*),
\ {\rm as}\ \Delta\to\infty.
\EQNY
%Finally, with $\Delta \to \IF$ we get
%\BQNY
%\lefteqn{\sum_{k=1}^{N_u}\sum_{l=-N_u}^{N_u}\pk*{\exists_{s \in E_{u,k}^1, t \in E_{u,l}^2}: W_1^*(s)>u, W_2^*(t)> au  }}\\
%\EQNY
To complete the proof, (3.14) needs to be shown to be asymptotically negligible, which is given below in Section 9.

%%%%%%%%%%%%%%%%%%%%%%%%%%%%%%%%%%%%%%%%%%%%%%%%%%%%%

\section{Proof of negligibility  of (3.14)}
For any $-N_u \le l \le N_u$
\BQNY
\lefteqn{\sum_{k=1}^{N_u-l}\pk*{\exists_{s \in E_{u,1}^1, t_1 \in E_{u,k+l}^2, t_2 \in E_{u,l}^2}:
		\begin{array}{ccc}
			W_1^*(s)>u \\
			W_2^*(t_1) > a u \\
			W_2^*(t_2) > au
		\end{array}
}}\\
&= & \sum_{k=1}^{N_u-l}\int_{\R}\phi(u+c_1-\frac{x}{u})\\
&&\times \pk*{\exists_{s \in E_{u,0}^1, t_1 \in E_{u,k+l}^2, t_2 \in E_{u,l}^2}:
	\begin{array}{ccc}
		W_1^*(s)>u \\
		W_2^*(t_1) > a u \\
		W_2^*(t_2) > au
	\end{array}
	\Bigg{|}W_1(1)=u+c_1-\frac{x}{u}}dx \\
& = & \sum_{k=1}^{N_u-l}\int_{\R}\phi(u+c_1-\frac{x}{u})\pk*{\exists_{s \in E_{u,0}^1}:W_1(s)-W_1(1)+c_1-c_1 s >\frac{x}{u} }\\
&&\times \pk*{\exists_{t_1 \in E_{u,k+l}^2, t_2 \in E_{u,l}^2}
	\begin{array}{ccc}
		W_2^*(t_1) > a u \\
		W_2^*(t_2) > au
	\end{array}
	\Bigg{|}W_1(1)=u+c_1-\frac{x}{u}}dx \\
& = & \sum_{k=1}^{N_u-l}\int_{\R}\phi(u+c_1-\frac{x}{u})\pk*{\exists_{s \in E_{u,0}^1}:W_1(s)-W_1(1)+c_1-c_1 s >\frac{x}{u} }\\
&& \times \pk*{\exists_{t_1 \in E_{u,k+l}^2, t_2 \in E_{u,l}^2}:X_{x,u}(t_1,t_2)>0}dx,
\EQNY
where $X_{x,u}(t_1,t_2)=(X_{1,x,u}(t_1),X_{2,x,u}(t_2))$ is a bivariate Gaussian process with \[\E {X_{x,u}(t_1,t_2)}=-\Bigg{(}\begin{array}{ccc}
-c_2t_1+\rho t_1 (c_1-\frac{x}{u}) \\
-c_2t_2+\rho t_2 (c_1-\frac{x}{u})
\end{array}\Bigg{)}+\Bigg{(}\begin{array}{ccc}
-(a-\rho t_1)u \\
-(a-\rho t_2)u
\end{array}\Bigg{)}\] and
\[\Sigma_{X_{x,u}(s,t)}=\left(\begin{array}{ccc}
t_1-\rho^2 t_1^2 & t_1-\rho^2 t_1 t_2 \\
t_1-\rho^2 t_1 t_2 & t_2-\rho^2 t_2^2
\end{array}\right).\]
Denote
$$S_0=\pk*{\exists_{t_1 \in E_{u,l}^2, t_2 \in E_{u,l}^2}:X_{x,u}(t_1,t_2)>0}, S_1=\sum_{k=2}^{N_u-l}\pk*{\exists_{t_1 \in E_{u,k+l}^2, t_2 \in E_{u,l}^2}:X_{x,u}(t_1,t_2)>0},$$
$$S_2=\pk*{\exists_{t_1 \in (1-\frac{(l+2)\Delta}{u^2},1-\frac{(l+1+\frac{1}{\sqrt{\Delta}})\Delta}{u^2}), t_2 \in E_{u,l}^2}:X_{x,u}(t_1,t_2)>0},$$
$$S_3=\pk*{\exists_{t_1 \in(1-\frac{(l+1+\frac{1}{\sqrt{\Delta}})\Delta}{u^2},1-\frac{(l+1)\Delta}{u^2}), t_2 \in E_{u,l}^2}:X_{x,u}(t_1,t_2)>0}.$$
Observe that for (3.14) to be negligible it is enough to show that, as $u \to \infty$
$$\frac{S_1+S_2+S_3}{S_0} \to 0.$$

Notice that for $X_{x,u}(t_1,t_2)=(X_{1,x,u}(t_1),X_{2,x,u}(t_2))$ we have
\BQNY
\pk*{\exists_{t_1 \in E_{u,k+l}^2, t_2 \in E_{u,l}^2}:X_{x,u}(t_1,t_2)>0} &\le& \pk*{\exists_{t_1 \in E_{u,k+l}^2, t_2 \in E_{u,l}^2}:
	X_{1,x,u}(t_1)+X_{2,x,u}(t_2)>0} \\
& \le & \pk*{\exists_{t_1 \in E_{u,k+l}^2, t_2 \in E_{u,l}^2}:
	\frac{X_{1,x,u}(t_1)+X_{2,x,u}(t_2)}{\sigma_{k,u}}>0},
\EQNY
where $$\sigma_{k,u}^2:=\max_{t_1 \in E_{u,k+l}^2, t_2 \in E_{u,l}^2} \eta_{u}^2(t_1,t_2) $$
and $\eta_{u}^2(t_1,t_2):=Var(X_{1,x,u}(t_1)+X_{2,x,u}(t_2)).$
Then for any $t_1 \in E_{u,k+l}^2, t_2 \in E_{u,l}^2$
$$\frac{\partial \eta_{u}^2(t_1,t_2)}{\partial t_1}=3-(2\rho^2+2t_2\rho^2), \frac{\partial \eta_{k,u}^2(t_1,t_2)}{\partial t_2}=1-(2\rho^2+2t_1\rho^2)$$
as $u \to \infty.$
\begin{enumerate}
	\item If $\frac{\partial \eta_{u}^2(t_1,t_2)}{\partial t_1}>0, \frac{\partial \eta_{u}^2(t_1,t_2)}{\partial t_2}>0,$ then $t_1^*=\frac{a}{\rho(2a\rho-1)}-\frac{(l+k)\Delta}{u^2}, t_2^*=\frac{a}{\rho(2a\rho-1)}-\frac{l\Delta}{u^2}.$ Consequently
	\BQN \label{Sig}
	\sigma_{k,u}^2&=&-\frac{1}{\rho(2a\rho-1)^2}\left(4a-4a^2\rho^2-\frac{1}{u^2}(f_k+g_l)+O(\frac{1}{u^4})\right),
	\EQN
	where $g_l=8al\Delta \rho^2-4l\Delta\rho.$
	\item If $\frac{\partial\ \eta_{u}^2(t_1,t_2)}{\partial t_1}>0, \frac{\partial \eta_{u}^2(t_1,t_2)}{\partial t_2}<0,$ then $t_1^*=\frac{a}{\rho(2a\rho-1)}-\frac{(l+k)\Delta}{u^2}, t_2^*=\frac{a}{\rho(2a\rho-1)}-\frac{(l+1)\Delta}{u^2}.$ Consequently \eqref{Sig} holds with
	$g_l=8al\Delta \rho^2-4l\Delta\rho-\Delta \rho+4a^2\Delta\rho^3.$
	\item If $\frac{\partial \eta_{u}^2(t_1,t_2)}{\partial t_1}<0, \frac{\partial \eta_{u}^2(t_1,t_2)}{\partial t_2}<0,$ then $t_1^*=\frac{a}{\rho(2a\rho-1)}-\frac{(l+k+1)\Delta}{u^2}, t_2^*=\frac{a}{\rho(2a\rho-1)}-\frac{l\Delta}{u^2}.$ Consequently \eqref{Sig} holds with
	$g_l=8al\Delta \rho^2-4l\Delta\rho-4\Delta \rho+8a\Delta \rho^2.$
\end{enumerate}
In all of the above scenarios $f_k=8ak\Delta \rho^2-3k\Delta\rho-4a^2k\Delta \rho^3.$ Notice that for $1>\frac{a}{\rho(2a\rho-1)}>0$ we have $f_k>0$ and $\rho<0.$ Denote $\mu_u:=\E{X_1(t_1^*)+X_2(t_2^*)}=2au+c_2t_1^*+c_2t_2^*-\rho (t_1^*+t_2^*)(u+c_1-\frac{x}{u})$. For all $i \in {1,2,3},$ using \cite{Pit96}[Thm 8.1], there exist constants $C, C_2>0$ such that
\BQN \label{s1IV}
S_1 &\le& \sum_{k=2}^{N_u-l} C\frac{\mu_u}{\sigma_{k,u}} e^{-\frac{\mu_u^2}{2\sigma_{k,u}^2}}\nonumber\\
&=& \sum_{k=2}^{N_u-l} C \frac{\mu_u}{\sigma_{k,u}}e^{-\frac{\mu_u^2\left(4a-4a^2\rho^2+\frac{1}{u^2}(f_k+g_l^{(i)})+O(\frac{1}{u^4})\right)}{-\frac{2}{\rho(2a\rho-1)^2}((4a-4a^2\rho^2)^2+O(\frac{1}{u^4}))}}\nonumber\\
&\le&C \frac{\mu_u}{\sigma_{0,u}}e^{-\frac{\mu_u^2(4-4\rho^2+\frac{g_l^{(i)}}{u^2}+O(\frac{1}{u^4}))}{-\frac{2}{\rho(2a\rho-1)^2}((4a-4a^2\rho^2)^2+O(\frac{1}{u^4}))}}\sum_{k=2}^{N_u-l} e^{-C_2 k( \Delta+O(\frac{1}{u^2}))}\nonumber\\
&\le&C \frac{\mu_u}{\sigma_{0,u}}e^{-\frac{\mu_u^2(4-4\rho^2+\frac{g_l^{(i)}}{u^2}+O(\frac{1}{u^4}))}{-\frac{2}{\rho(2a\rho-1)^2}((4a-4a^2\rho^2)^2+O(\frac{1}{u^4}))}}\frac{e^{-C_2\Delta}}{e^{C_2\Delta}-1}.
\EQN
Similarly we get that
\BQN \label{s2IV}
S_2 &\le&C \frac{\mu_u}{\sigma_{0,u}}e^{-\frac{\mu_u^2(4-4\rho^2+\frac{g_l^{(i)}}{u^2}+O(\frac{1}{u^4}))}{-\frac{1}{\rho(2a\rho-1)^2}((4a-4a^2\rho^2)^2+O(\frac{1}{u^4}))}}e^{-C_2\sqrt{\Delta}}
\EQN
and
\BQN \label{s3IV}
S_3 &\le & \pk*{\exists_{t_1 \in(1-\frac{(l+1+\frac{1}{\sqrt{\Delta}})\Delta}{u^2},1-\frac{(l+1)\Delta}{u^2})}:X_1(t_1)>0} \nonumber \\
& \le & \frac{\sqrt{\Delta}}{\Delta}\pk*{\exists_{t_1 \in(1-\frac{(l+2)\Delta}{u^2},1-\frac{(l+1)\Delta}{u^2})}:X_1(t_1)>0}.
\EQN
Using (\ref{s1IV}), (\ref{s2IV}) and (\ref{s3IV}) we have that for some $C>0$
$$\frac{S_1+S_2+S_3}{S_0} \le \frac{C}{\sqrt{\Delta}}+e^{-C \sqrt{\Delta}}+\frac{e^{-C\Delta}}{e^{C\Delta}-1} \to 0, \Delta \to \infty.  $$
Hence the proof follows. \QED

%%%%%%%%%%%%%%%%%%%%%%%%%%%%%%%%%%%%%%%%%%%%%%%%%%%%%

\section{Proof of negligibility of (3.15)}
By taking $ \vk b= \Sigma^{-1} _{1,t^* }  (1,a)>(0,0)$,
asymptotics given in the proof of case (v) of Theorem 2.2 imply that
\BQNY
\lim_{u \to \IF}\frac{1}{u^2}\log \pk*{ \exists_{(s,t) \in F_{i,u}}: W_1^*(s) > u, W_2^*(t)>u}
=-\frac{1}{2V_1},
%Var\left(\frac{b_1W_1^*(1)+b_2W_2^*(t^*)}{2(b_1+b_2)} + \frac{b_1W_1^*(t^*)+b_2W_2^*(1)}{2(b_1+b_2)}\right)}.
\EQNY
for $i=1,2$, where
$V_1:=Var\left(\frac{b_1W_1^*(1)+b_2W_2^*(t^*)}{b_1+b_2}\right)= Var\left(\frac{b_2W_1^*(t^*)+b_1W_2^*(1)}{b_1+b_2}\right)$.

Moreover
\BQN
\pk*{ \exists_{(s,t) \in F_{1,u},(s',t') \in F_{2,u}}: W_1^*(s) > u, W_2^*(t)>u,W_1^*(s') > u, W_2^*(t')>u}\nonumber\\
&&\hspace*{-12cm}\le \pk*{\exists_{(s,t) \in F_{1,u},(s',t') \in F_{2,u}}: \frac{b_1W_1^*(s)+b_2W_2^*(t)}{2(b_1+b_2)} + \frac{b_2W_1^*(s')+b_1W_2^*(t')}{2(b_1+b_2)}>u}.
\label{proc_sup}
\EQN
Since
$$\lim_{u \to \IF}F_{1,u}=\{(1,t^*)\}, \quad \lim_{u \to \IF}F_{2,u}=\{(t^*,1)\}$$
and  variance function of process under supremum in (\ref{proc_sup}) is continuous,
then using Borell-TIS inequality (see e.g., \cite{Pit96}) we get
\BQNY
\lim_{u \to \IF}\frac{1}{u^2}\log\pk*{\exists_{(s,t) \in F_{1,u},(s',t') \in F_{2,u}}: \frac{b_1W_1^*(s)+b_2W_2^*(t)}{2(b_1+b_2)}
	+ \frac{b_2W_1^*(s')+b_1W_2^*(t')}{2(b_1+b_2)}>u}
&\le&-\frac{1}{2V_2},
%Var\left(\frac{b_1W_1^*(1)+b_2W_2^*(t^*)}{2(b_1+b_2)} + \frac{b_1W_1^*(t^*)+b_2W_2^*(1)}{2(b_1+b_2)}\right)}.
\EQNY
where
$$V_2:=Var\left(\frac{b_1W_1^*(1)+b_2W_2^*(t^*)}{2(b_1+b_2)} + \frac{b_2W_1^*(t^*)+b_1W_2^*(1)}{2(b_1+b_2)}\right).$$

Since $t^* <1$
%Next define
%$$V_1:=Var\left(\frac{b_1W_1^*(1)+b_2W_2^*(t^*)}{b_1+b_2}\right),$$
%$$V_2:=Var\left(\frac{b_1W_1^*(1)+b_2W_2^*(t^*)}{2(b_1+b_2)} + \frac{b_1W_1^*(t^*)+b_2W_2^*(1)}{2(b_1+b_2)}\right).$$
%Observe that
%\BQNY
%\lim_{u \to \IF}\frac{1}{u^2}\log\left(\frac{\pk*{ \exists_{(s,t) \in F_{1,u}}\exists_{(s',t') \in F_{2,u}} W_1^*(s) > u, W_2^*(t)>u,W_1^*(s') > u, W_2^*(t')>u}}{\pk*{ \exists_{(s,t) \in F_{1,u}} W_1^*(s) > u, W_2^*(t)>u}}\right)
%&=&-\frac{1}{2}(\frac{1}{V_2}-\frac{1}{V_1})\\
%&=&\frac{1}{2}\frac{V_2-V_1}{V_1V_2}
%\EQNY
%Hence for \eqref{quatro} to be negligible it suffices to show that $V_1 \ge V_2.$ Notice that because of symmetry
%$$
%V_2=\frac{1}{2}Var\left(\frac{b_1W_1^*(1)+b_2W_2^*(t^*)}{b_1+b_2}\right)+\frac{2b_1b_2t^*+b_1^2\rho+b_2^2\rho t^*}{2(b_1+b_2)^2}.
%$$
%Therefore since $t^* <1,$ then
\BQNY
V_1-V_2 %&=&\frac{1}{2}Var\left(\frac{b_1W_1^*(1)+b_2W_2^*(t^*)}{b_1+b_2}\right)-\frac{2b_1b_2t^*+b_1^2\rho+b_2^2\rho t^*}{2(b_1+b_2)^2}\\
&=&\frac{b_1^2(1-\rho)+b_2^2t^*(1-\rho)-2b_1b_2t^*(1-\rho)}{2(b_1+b_2)^2}\\
&>&(1-\rho)\frac{b_1^2+b_2^2t^*-2b_1b_2\sqrt{t^*}}{2(b_1+b_2)^2}\\
&=&(1-\rho)\frac{(b_1-b_2\sqrt{t^*})^2}{2(b_1+b_2)^2} \ge 0.
\EQNY
Hence (3.15) is asymptotically negligible as $u\to\infty$.

\bibliographystyle{ieeetr}

\bibliography{queue2d}
\end{document}